\documentclass[11pt]{article}
\usepackage{graphicx}
\usepackage{amssymb}
\usepackage{amsmath}
\usepackage{amsthm}
\newtheorem{definition}{Definition}[section]
\newtheorem{lemma}{Lemma}[section]
\newtheorem{theorem}{Theorem}[section]
\newtheorem{proposition}{Proposition}[section]
\newtheorem{corollary}{Corollary}[section]
\newtheorem{oss}{Remark}[section]

{\end{description}}
\newenvironment{example}{\begin{description} \item[\bf Example:] }%
{\end{description}}

{\hfill $\bullet$ \\}

\newcommand{\adots}{\cdot^{\textstyle\cdot^{\textstyle\cdot}}}

\newcommand{\be}{\begin{equation}}
\newcommand{\ee}{\end{equation}}

\newcommand{\norma}[1]{\left\Vert #1 \right\Vert}

\newcommand{\Test}{{{\cal C}_0}}
\newcommand{\diag}{\begin{smallmatrix}\vspace{-0.5ex}\textrm{\normalsize diag}\\\vspace{-0.8ex}j=1,\ldots,n\end{smallmatrix}}

\begin{document}
\title{Spectral features and asymptotic properties for $\alpha$-circulants and $\alpha$-Toeplitz sequences:
theoretical results and examples}
\author{Eric Ngondiep, Stefano Serra-Capizzano, and Debora Sesana
       \thanks{Dipartimento di Fisica e
       Matematica, Universit\`a dell'Insubria,
       Via Valleggio 11, 22100 Como (ITALY).
       Email:\ \{eric.ngondiep,stefano.serrac,debora.sesana\}@uninsubria.it; serra@mail.dm.unipi.it}
       }
\maketitle
\date{}

\begin{abstract}
For a given nonnegative integer $\alpha$, a matrix $A_n$ of size $n$
is called $\alpha$-Toeplitz if its entries obey the rule
$A_n=\left[a_{r-\alpha s}\right]_{r,s=0}^{n-1}$. Analogously, a
matrix $A_n$ again of size $n$ is called $\alpha$-circulant if
$A_n=\left[a_{(r-\alpha s)\ {\rm mod}\, n}\right]_{r,s=0}^{n-1}$.
Such kind of matrices arises in wavelet analysis, subdivision
algorithms and more generally when dealing with multigrid/multilevel
methods for structured matrices and approximations of boundary value
problems. In this paper we study the singular values of
$\alpha$-circulants and we provide an asymptotic analysis of the
distribution results for the singular values of $\alpha$-Toeplitz
sequences in the case where $\{a_k\}$ can be interpreted as the sequence of
Fourier coefficients of an integrable function $f$ over the domain
$(-\pi,\pi)$. Some generalizations to the block, multilevel case,
amounting to choose $f$ multivariate and matrix valued, are briefly
considered.
\end{abstract}
\\
{\bf AMS SC:} 65F10, 15A18.

\section{Introduction} \label{sec1}

A matrix $A_n$ of size $n$ is called $\alpha$-Toeplitz if its
entries obey the rule $A_n=\left[a_{r-\alpha
s}\right]_{r,s=0}^{n-1}$, where $\alpha$ is a nonnegative integer.
As an example, if $n=5$ and $\alpha=3$ then
\begin{eqnarray*}
  A_n\equiv T_{n,\alpha}=\left[\begin{array}{ccccc}
    a_{0} & a_{-3} & a_{-6} & a_{-9} & a_{-12} \\
    a_{1} & a_{-2} & a_{-5} & a_{-8} & a_{-11} \\
    a_{2} & a_{-1} & a_{-4} & a_{-7} & a_{-10} \\
    a_{3} & a_{0}  & a_{-3} & a_{-6} & a_{-9} \\
    a_{4} & a_{1}  & a_{-2} & a_{-5} & a_{-8}
  \end{array}\right].
\end{eqnarray*}

Along the same lines, a matrix $A_n$ of size $n$ is called
$\alpha$-circulant if $A_n=\left[a_{(r-\alpha s)\ {\rm mod}\,
n}\right]_{r,s=0}^{n-1}$. For instance if $n=5$ and $\alpha=3$ then
we have
\begin{eqnarray*}
  A_n\equiv C_{n,\alpha}=\left[\begin{array}{ccccc}
    a_{0} & a_{2} & a_{4} & a_{1} & a_{3} \\
    a_{1} & a_{3} & a_{0} & a_{2} & a_{4} \\
    a_{2} & a_{4} & a_{1} & a_{3} & a_{0} \\
    a_{3} & a_{0} & a_{2} & a_{4} & a_{1} \\
    a_{4} & a_{1} & a_{3} & a_{0} & a_{2}
  \end{array}\right].
\end{eqnarray*}

Such kind of matrices arises in wavelet analysis \cite{wave} and
subdivision algorithms or, equivalently, in the associated
refinement equations, see \cite{subd} and references therein.
Furthermore, it is interesting to remind that Gilbert Strang
\cite{strang} has shown rich connections between dilation
equations in the wavelets context and multigrid methods
\cite{Hack,Trot}, when constructing the restriction/prolongation
operators \cite{FS2,ADS} with various boundary conditions. It is
worth noticing that the use of different boundary conditions is
quite natural when dealing with signal/image restoration problems or
differential equations, see \cite{model-tau,Sun}.

In this paper we address the problem of characterizing the singular
values of $\alpha$-circulants and of providing an asymptotic
analysis of the distribution results for the singular values of
$\alpha$-Toeplitz sequences, in the case where the sequence of
values $\{a_k\}$, defining the entries of the matrices, can be interpreted as the sequence of
Fourier coefficients of an integrable function $f$ over the domain
$(-\pi,\pi)$. As a byproduct, we will show interesting relations
with the analysis of convergence of multigrid methods given, e.g.,
in \cite{mcirco,ADS}. Finally we generalize the analysis to the
block, multilevel case, amounting to choose the symbol $f$
multivariate, i.e., defined on the set $(-\pi,\pi)^d$ for some
$d>1$, and matrix valued, i.e., such that $f(x)$ is a matrix of
given size $p\times q$.

The paper is organized as follows. In Section \ref{sec:tools} we
report useful definitions, well-known results in the standard case
of circulants and Toeplitz that is when $\alpha=1$ (or $\alpha=e$, $e=(1,\ldots,1)$, in the multilevel setting), and a
preliminary analysis of some special cases. Section \ref{sec:circ}
deals with the singular value analysis of $\alpha$-circulants while
in Section \ref{sec:toep} we treat the $\alpha$-Toeplitz case in an
asymptotic setting, and more precisely in the sense of the Weyl
spectral distributions. Section \ref{sec:multigrid} is devoted to
sketch useful connections with multigrid methods, while in Section
\ref{sec:gen} we report the generalization of the results when we
deal the multilevel block case. Section \ref{sec:fin} is aimed to
draw conclusions and to indicate future lines of research.

\section{General definitions and tools}\label{sec:tools}

For any $n\times n$ matrix $A$ with eigenvalues $\lambda_j(A)$,
$j=1,\ldots,n$, and for any $m\times n$ matrix $B$ with singular
values $\sigma_j(B)$, $j=1,\ldots,l$, $l=\min\{m,n\}$, we set
\[ {\rm Eig}(A)= \{\lambda_1(A), \lambda_2(A),\ldots,\lambda_n(A)\},
\ \ \ \ {\rm Sgval}(B)= \{\sigma_1(B),
\sigma_2(B),\ldots,\sigma_l(B)\}.\]

The matrix $B^{*}B$ is positive semidefinite, since
$x^{*}(B^{*}B)x=\|Bx\|_{2}^{2}\geq 0$ for all $x\in \mathbb{C}^{n}$,
with $^*$ denoting the transpose conjugate operator. Moreover, it is
clear that the eigenvalues $\lambda_{1}(B^*B)\geq
\lambda_{2}(B^*B)\geq\cdots\geq \lambda_{n}(B^*B)\geq 0$ are
nonnegative and can therefore be written in the form
\begin{equation}\label{1}
\lambda_{j}(B^*B)=\sigma_{j}^{2},
\end{equation}
with $\sigma_{j}\geq0$, $j=1,\ldots,n$. The numbers
$\sigma_{1}\geq\sigma_{2}\geq\cdots\geq\sigma_{l}\geq0$,
$l=\min\{m,n\}$, are called singular values of $B$, i.e.,
$\sigma_{j}=\sigma_{j}(B)$ and if $n>l$ then $\lambda_{j}(B^*B)=0$,
$j=l+1,\ldots,n$. A more general statement is contained in the
singular value decomposition theorem (see e.g. \cite{GV}).
% \begin{oss}
%The smallest singular value $\sigma_{n}$ of a square matrix $B$
%(when $m=n$) gives the distance of $B$ to the ``nearest'' singular
%matrix. \end{oss}

\begin{theorem}\label{2}
Let $B$ be an arbitrary (complex) $m\times n$ matrix. Then:
\begin{description}
    \item (a) There exists a unitary $m\times m$ matrix $U$ and a
    unitary $n\times n$ matrix $V$ such that $U^{*}BV=\Sigma$ is an
    $m\times n$ ``diagonal matrix'' of the following form:
    \begin{equation*}
        \Sigma=\begin{bmatrix}
                 D & 0 \\
                 0 & 0 \\
               \end{bmatrix}, \text{\,\,}
               D:={\rm diag}(\sigma_{1},\ldots,\sigma_{r}),\text{\,\,\,}
               \sigma_{1}\geq\sigma_{2}\geq\cdots\geq\sigma_{r}>0.
    \end{equation*}

    Here $\sigma_{1},\ldots,\sigma_{r}$ are the nonvanishing singular
    values of $B$, and $r$ is the rank of $B.$
    \item(b) $\;\,$The nonvanishing singular values of $B^{*}$ are also
    precisely the number $\sigma_{1},\ldots,\sigma_{r}.$\\
    The decomposition $B=U\Sigma V^{*}$ is called
    ``the singular value decomposition of $B$''.
\end{description}
\end{theorem}
% \begin{theorem}\label{BAB} Let $A$ be an arbitrary (complex) $n\times n$ matrix, and $B$ a
%unitary $n\times n$ matrix, then $A$ and $B^{*}AB$ have the same eigenvalues. \end{theorem}
% \begin{oss} The tree theorems are proved by \textbf{Schur}. \end{oss}

For any function $F$ defined on ${\mathbb R}^+_0$ and for any
$m\times n$ matrix $A$, the symbol $\Sigma_{\sigma}(F,A)$ stands for
the mean
\begin{eqnarray}\label{sigmaFA}
\Sigma_{\sigma}(F,A):= {\frac 1 {\min\{n,m\}}
\sum_{j=1}^{\min\{n,m\}} F\left(\sigma_j(A)\right)}={\frac 1
{\min\{n,m\}} \sum_{\sigma\in {\rm Sgval}(A)} F(\sigma)}.
\end{eqnarray}

Throughout this paper we speak also of {\em matrix sequences} as
sequences $\{A_k\}$ where $A_k$ is an $n(k)\times m(k)$ matrix with
$\min\{n(k),m(k)\}\rightarrow \infty$ as $k\rightarrow \infty$. When
$n(k)=m(k)$ that is all the involved matrices are square, and this
will occur often in the paper, we will not need the extra parameter
$k$ and we will consider simply matrix sequences of the form
$\{A_n\}$.

Concerning the case of matrix-sequences an important notion is that
of spectral distribution in the eigenvalue or singular value sense,
linking the collective behavior of the eigenvalues or singular
values of all the matrices in the sequence to a given function (or
to a measure). The notion goes back to Weyl and has been
investigated by many authors in the Toeplitz and Locally Toeplitz
context (see the book by B\"ottcher and Silbermann \cite{BS} where
many classical results by the authors, Szeg\"o, Avram, Parter, Widom
Tyrtyshnikov, and many other can be found, and more recent results
in \cite{jacobi-GOL,ku-ser,zabroda,tyrtL1,Tillinota,tillicomplex}).
Here we report the definition of spectral distribution only in the
singular value sense since our analysis is devoted to singular
values. The case of eigenvalues will be the subject of future
investigations.

\begin{definition}\label{def-distribution}
{\rm Let $\mathcal C_0({\mathbb R}^+_0)$ be the set of continuous
functions with bounded support defined over the nonnegative real
numbers, $d$ a positive integer, and $\theta$ a complex-valued
measurable function defined on a set $G\subset\mathbb R^d$ of finite
and positive Lebesgue measure $\mu(G)$. Here $G$ will be often equal
to $(-\pi,\pi)^d$ so that $e^{i\overline{G}}={\mathbb T}^d$
with ${\mathbb T}$ denoting the complex unit circle. A matrix
sequence $\{A_k\}$ is said to be {\em distributed $($in the sense of
the singular values$)$ as the pair $(\theta,G)$,} or to {\em have
the distribution function $\theta$} ($\{A_k\}\sim_{\sigma}(\theta,G)$), if, $\forall F\in \mathcal
C_0({\mathbb R}^+_0)$, the following limit relation holds
\begin{equation}\label{distribution:sv-eig}
\lim_{k\rightarrow
\infty}\Sigma_{\sigma}(F,A_k)=\frac1{\mu(G)}\,\int_G F(|\theta(t)|)\,
dt,\qquad t=(t_{1},\ldots,t_{d}).
\end{equation}

When considering $\theta$ taking values in  ${\cal {M}}_{pq}$, where
${\cal {M}}_{pq}$ is the space of $p \times q$ matrices with complex
entries and a function is considered to be measurable if and only if
the component functions are, we say that $\{A_k\}\sim_{\sigma}
(\theta,G)$ when for every $F\in \mathcal C_0({\mathbb R}^+_0)$ we
have
\[
\lim_{k\rightarrow \infty}\Sigma_{\sigma}(F,A_k)= \frac {1}
{\mu(G)}\,\int_G \frac {\sum_{j=1}^{\min\{p,q\}}
\left(F(\sigma_j(\theta(t)))\right)} {\min\{p,q\}}\, dt,\qquad t=(t_{1},\ldots,t_{d}),
\]
with $\sigma_j(\theta(t))=\sqrt{\lambda_j(\theta(t)\theta^*(t))}=\lambda_j(\sqrt{\theta(t)\theta^*(t)})$.
Finally we say that two sequences $\{A_k\}$ and $\{B_k\}$ are {\em
equally distributed} in the sense of singular values ($\sigma$) if,
$\forall F\in \mathcal C_0({\mathbb R}^+_0)$, we have
\[
\lim_{k\rightarrow\infty}[\Sigma_{\sigma}(F,B_k)-\Sigma_{\sigma}(F,A_k)]=0.
\]
}
\end{definition}

Here we are interested in explicit formulae for the singular values
of $\alpha$-circulants and in distribution results for
$\alpha$-Toeplitz sequences. In the latter case, following what is
known in the standard case of $\alpha=1$ (or $\alpha=e$ in the multilevel setting), we need to link the
coefficients of the $\alpha$-Toeplitz sequence to a certain symbol.

Let $f$ be a Lebesgue integrable function defined on $(-\pi,\pi)^d$
and taking values in ${\cal {M}}_{pq}$, for given positive integers
$p$ and $q$. Then, for $d$-indices $r=(r_1,\ldots,r_d),
j=(j_1,\ldots,j_d), n=(n_1,\ldots,n_d)$, $e=(1,\ldots,1)$, $\underline{0}=(0,\ldots,0)$, the Toeplitz matrix
$T_n(f)$ of size $p\hat n\times q\hat n$, $\hat n=n_1\cdot n_2\cdots
n_d$, is defined as follows
\begin{eqnarray*}%\label{toepmat}
T_{n}(f)=[\tilde{f}_{r-j}]_{r,j=\underline{0}}^{n-e},
\end{eqnarray*}
where $\tilde{f}_{k}$ are the Fourier coefficients of $f$ defined by
equation
\begin{equation}\label{defcoeff}
\tilde{f}_j=\tilde{f}_{(j_1, \ldots, j_d)} (f) = \frac 1 {(2\pi)^d}
\int_{{[-\pi,\pi]}^d} f(t_1,\ldots,t_d)e^{-i(j_1t_1 + \cdots +
j_dt_d)}\, dt_1 \cdots dt_d,\quad \quad i^2=-1,
\end{equation}
for integers $j_{\ell}$ such that $-\infty < j_{\ell} < \infty$ for
$1 \le \ell \le d$. Since $f$ is a matrix-valued function of $d$
variables whose component functions are all integrable, then the
$(j_1, \ldots, j_d)$-th Fourier coefficient is considered to be the
matrix whose $(u,v)$-th entry is the $(j_1, \ldots, j_d)$-th Fourier
coefficient of the function $(f(t_1,\ldots,t_d))_{u,v}$.

According to this multi-index block notation we can define general
multi-level block $\alpha$-Toeplitz and $\alpha$-circulants. Of
course, in this multidimensional setting, $\alpha$ denotes a
$d$-dimensional vector of nonnegative integers that is
$\alpha=(\alpha_1,\ldots,\alpha_d)$. In that case
$A_n=\left[a_{r-\alpha \circ s}\right]_{r,s=\underline{0}}^{n-e}$ where the
$\circ$ operation is the componentwise Hadamard product between
vectors or matrices of the same size. A matrix $A_n$ of size $p\hat
n\times q\hat n$ is called $\alpha$-circulant if
$A_n=\left[a_{(r-\alpha \circ s)\ {\rm mod}\,
n}\right]_{r,s=\underline{0}}^{n-e}$, where
  \begin{eqnarray*}
    (r-\alpha\circ s)\textrm{ mod $n$}=\left((r_{1}-\alpha_{1}s_{1})\textrm{ mod $n_{1}$},(r_{2}-\alpha_{2}s_{2})\textrm{ mod $n_{2}$},\ldots,(r_{d}-\alpha_{d}s_{d})\textrm{ mod $n_{d}$}\right).
  \end{eqnarray*}

\subsection{The extremal cases where $\alpha=\underline{0}$ or $\alpha=e$, and the intermediate cases}

We consider a $d$-level setting and we analyze in detail the case where $\underline{0}\le
\alpha\le e$ and with $\le$ denoting the componentwise partial
ordering between real vectors. When $\alpha$ has at least a zero component, the analysis can be reduced to the positive one as studied in Subsection \ref{nonnegative-vs-positive}.

\subsubsection{$\alpha=e$}

In the literature the only case deeply studied is the case of
$\alpha=e$ (standard shift in every level). Here for multilevel
block circulants $A_{n}=[a_{(r-\alpha\circ s)\textrm{ mod $n$}}]_{r,s=\underline{0}}^{n-e}$ the singular values are given by those of
\begin{eqnarray*}
  \sigma_{k}(A_{n})=\sum_{j=\underline{0}}^{n-e} a_j e^{i2\pi(j_1 k_1/n_1  + \cdots
+ j_d k_d/n_d)},\qquad k=(k_{1},\ldots,k_{d}),
\end{eqnarray*}
for any $k_{\ell}$ such that $0\le k_{\ell}\le n_{\ell}-1$,
$\ell=1,\ldots,d$. Of course when the coefficients $a_j$ comes from
the Fourier coefficients of a given Lebesgue integrable function
$f$, i.e. $\tilde{f}_j=a_{j\ {\rm mod}\, n}$, $j=-n/2,\ldots,n/2$
($n/2=(n_{1}/2,n_{2}/2,\ldots,n_{d}/2)$), the singular values are
those of $n/2$-th Fourier sum of $f$ evaluated at the grid points
\[
2\pi k/n=2\pi\left(k_1/n_1,\ldots,k_d/n_d\right),
\]
$0\le k_j\le n_j-1$, $j=1,\ldots,d$. Moreover the explicit Schur
decomposition is known. For $d=p=q=1$ any standard circulant matrix
can be written in the form
\begin{eqnarray}\label{iV}
  A_{n}\equiv C_{n} =F_{n}D_{n}F_{n}^{\ast},
\end{eqnarray}
where
\begin{eqnarray}\label{V}
 \notag F_{n} &=&\frac{1}{\sqrt{n}}\left[e^{-\frac{2\pi ijk}{n}}\right]_{j,k=0}^{n-1},\text{\,\,Fourier matrix,} \\
  D_{n} &=& {\rm diag}(\sqrt{n}F_{n}^{\ast}\underline{a}),  \\
 \notag \underline{a} &=&\left[a_{0},a_{1},\ldots,a_{n-1}\right]^{T}, \text{\,\,\,
  first column of the matrix $A_{n}$}.
\end{eqnarray}

Of course for general $d,p,q$ the formula generalizes as
\[
 A_{n} =(F_{n}\otimes I_p) D_{n}(F_{n}^{\ast}\otimes I_q),
\]
with $F_n=F_{n_1}\otimes F_{n_2} \otimes \cdots \otimes F_{n_d}$
$D_{n} = {\rm diag}(\sqrt{\hat n}(F_{n}^{\ast}\otimes
I_p)\underline{a})$, where $\hat{n}=n_{1}\cdot n_{2}\cdots n_{d}$ and $\underline{a}$ being the first ``column'' of
$A_n$ whose entries $a_j$, $j=(j_1,\ldots,j_d)$, ordered
lexicographically, are blocks of size $p\times q$.

For multilevel block Toeplitz sequences $\{T_{n}(f)\}$ generated by an integrable $d$ variate and
matrix valued symbol $f$ the singular values are not explicitly
known but we know the distribution in the sense of Definition
\ref{def-distribution}; see \cite{Tillinota}. More precisely we have
\begin{equation}\label{szego-tyrty}
  \{T_n(f)\}\sim_{\sigma}
(f,Q^d), \quad \ \ Q=(-\pi,\pi).
\end{equation}

\subsubsection{$\alpha=\underline{0}$}\label{alphazero}

The other extreme is represented by the case where $\alpha$ is the
zero vector. Here the multilevel block $\alpha$-circulant and
$\alpha$-Toeplitz coincide when $\alpha=\underline{0}$ and are both given by
\begin{eqnarray*}
  A_{n}=[a_{(r-\underline{0}\circ s)\textrm{ mod $n$}}]_{r,s=\underline{0}}^{n-e}=[a_{r\textrm{ mod $n$}}]_{r,s=\underline{0}}^{n-e}=
    [a_{r}]_{r,s=\underline{0}}^{n-e}=\left[\begin{array}{ccc}
    a_{\underline{0}} & \cdots & a_{\underline{0}}\\
    \vdots &  & \vdots \\
    a_{n-e} & \cdots & a_{n-e}
    \end{array}\right].
\end{eqnarray*}

A simple computation shows that all the singular values are zero
except for few of them given by $\sqrt{\hat{n}}\sigma$, where $\hat{n}=n_{1}\cdot n_{2}\cdots n_{d}$ and
$\sigma$ is any singular value of the matrix $(\sum_{j=\underline{0}}^{n-e}
a_j^*a_j)^{1/2}$. Of course in the scalar case where $p=q=1$ the
choice of $\sigma$ is unique and by the above formula it coincides
with the Euclidean norm of the first column $\underline{a}$ of
the original matrix. In that case it is evident that
\begin{eqnarray*}
  \{A_n\}\sim_{\sigma}
(0,G),
\end{eqnarray*}
for any domain $G$ satisfying the requirements of Definition
\ref{def-distribution}.

\subsubsection{When some of the entries of $\alpha$ vanish}\label{nonnegative-vs-positive}

The content of this subsection reduces to the following remark: the
case of a nonnegative $\alpha$ can be reduced to the case of a
positive vector so that we are motivated to treat in detail the
latter in the next section. Let $\alpha$ be a $d$-dimensional vector
of nonnegative integers and let ${\cal{N}}\subset \{1,\ldots,d\}$ be
the set of indices such that $j\in \cal N$ if and only if
$\alpha_j=0$. Assume that $\cal N$ is nonempty, let $t\ge 1$ be its
cardinality and $d^+=d-t$. Then a simple calculation shows that the
singular values of the corresponding $\alpha$-circulant matrix $A_{n}=[a_{(r-\alpha\circ s)\textrm{ mod $n$}}]_{r,s=\underline{0}}^{n-e}$ are
zero except for few of them given by $\sqrt{\hat n[0]} \sigma$ where
\[
\hat n[0]=\prod_{j\in \cal{N}} n_j,\quad \quad
n[0]=(n_{j_1},\ldots,n_{j_t}),\quad {\cal{N}}=\{j_1,\ldots,j_t\},
\]
and $\sigma$ is any singular value of the matrix
\begin{equation}\label{eq-2-1-3}
\left(\sum_{j=\underline{0}}^{n[0]-e}C_j^*C_j\right)^{1/2}.
\end{equation}

Here $C_j$ is a $d^+$-level $\alpha^+$-circulant matrix with
$\alpha^+=(\alpha_{k_1},\ldots,\alpha_{k_{d^+}})$ and of partial
sizes $n[>0]=(n_{k_1},\ldots,n_{k_{d^+}})$,
${\cal{N}}^C=\{k_1,\ldots,k_{d^+}\}$, and whose expression is
\[
C_j=\left[a_{(r-\alpha \circ s)\ {\rm mod}\,
n}\right]_{r',s'=\underline{0}}^{n[>0]-e},
\]
where $(r-\alpha\circ s)_k=j_k$ for $\alpha_k=0$ and $r'_i=r_{k_i}$,
$s'_i=s_{k_i}$, $i=1,\ldots,d^+$. Taking into account the above
notation, for the $\alpha$-Toeplitz $A_{n}=[a_{r-\alpha\circ s}]_{r,s=\underline{0}}^{n-e}$ the same computation shows that
all the singular values are zero except for few of them given by
$\sqrt{\hat n[0]} \sigma$ where $\sigma$ is any singular value of the
matrix
\begin{equation}\label{eq-2-1-3-bis}
\left(\sum_{j=\underline{0}}^{n[0]-e}T_j^*T_j\right)^{1/2}.
\end{equation}

Here $T_j$ is a $d^+$-level $\alpha^{+}$-Toeplitz matrix with
$\alpha^+=(\alpha_{k_1},\ldots,\alpha_{k_{d^+}})$ and of partial
sizes $n[>0]=(n_{k_1},\ldots,n_{k_{d^+}})$,
${\cal{N}}^C=\{k_1,\ldots,k_{d^+}\}$, and whose expression is
\[
T_j=\left[a_{(r-\alpha \circ s)}\right]_{r',s'=\underline{0}}^{n[>0]-e},
\]
where $(r-\alpha\circ s)_k=j_k$ for $\alpha_k=0$ and $r'_i=r_{k_i}$,
$s'_i=s_{k_i}$, $i=1,\ldots,d^+$. Also in this case, since most of
the singular values are identically zero, we infer that
\begin{eqnarray*}
  \{A_n\}\sim_{\sigma}
(0,G),
\end{eqnarray*}
for any domain $G$ satisfying the requirements of Definition
\ref{def-distribution}.

\section{Singular values of $\alpha$-circulant matrices}\label{sec:circ}

Of course the aim of this paper is to give the general picture for
any nonnegative vector $\alpha$. Since the notations can become
quite heavy, for the sake of simplicity, we start with the case
$d=p=q=1$. Several generalizations, including also the
degenerate case in which $\alpha$ has some zero entries is treated
in Section \ref{sec:gen} via the observations in Subsection
\ref{nonnegative-vs-positive}, which imply that the general analysis
can be reduced to the case where all the entries of $\alpha$ are
positive, that is $\alpha_j>0$, $j=1,\ldots,d$.

In the following, we denote by $(n,\alpha)$ the greater common
divisor of $n$ and $\alpha$. i.e., $(n,\alpha)=\gcd(n,\alpha)$, by
$n_{\alpha}=\frac{n}{(n,\alpha)}$, by
$\check{\alpha}=\frac{\alpha}{(n,\alpha)}$, and by $I_t$ the
identity matrix of order $t$.

If we denote by $C_{n}$ the classical circulant matrix (i.e. with
$\alpha=1$) and by $C_{n,\alpha}$ the $\alpha$-circulant matrix
generated by its elements,  for generic $n$ and $\alpha$ one
verifies immediately that
\begin{eqnarray}\label{0}
  C_{n,\alpha} = C_{n}Z_{n,\alpha},
\end{eqnarray}
where
\begin{eqnarray}\label{i}
Z_{n,\alpha}=\left[\delta_{r-\alpha s}\right]_{r,s=0}^{n-1},\qquad\delta_{k}=
\left\{\begin{array}{cl} 1 & \textrm{if $k\equiv 0\textrm{ (mod $n$)}$,}\\0 &
\textrm{otherwise.}\end{array} \right.
\end{eqnarray}

\begin{lemma}
Let $n$ be any integer greater than $2$ then
\begin{eqnarray}\label{18}
  Z_{n,\alpha}=\underbrace{\left[\widetilde{Z}_{n,\alpha}|\widetilde{Z}_{n,\alpha}|\cdots|\widetilde{Z}_{n,\alpha}\right]}_{(n,\alpha)\textrm{ times}},
\end{eqnarray}
where $Z_{n,\alpha}$ is the matrix defined in $(\ref{i})$ and $\widetilde{Z}_{n,\alpha}\in \mathbb{C}^{n\times n_{\alpha}}$ is the
submatrix of $Z_{n,\alpha}$ obtained by considering only its first $n_{\alpha}$ columns, that is
\begin{eqnarray}\label{z-n-alfa}
  \widetilde{Z}_{n,\alpha} = Z_{n,\alpha}
  \left[\begin{array}{c}
  I_{n_{\alpha}} \\ 0 \end{array}\right].
\end{eqnarray}
\end{lemma}

\begin{proof} Setting $\widetilde{Z}_{n,\alpha}^{(0)}=\widetilde{Z}_{n,\alpha}$
and denoting by
$\widetilde{Z}_{n,\alpha}^{(j)}\in\mathbb{C}^{n\times n_{\alpha}}$
the $(j+1)$-th block-column of the matrix $Z_{n,\alpha}$ for
$j=0,\ldots,(n,\alpha)-1,$ we find
\begin{eqnarray*}
Z_{n,\alpha}=\left[\underbrace{\widetilde{Z}_{n,\alpha}^{(0)}}_{n\times n_{\alpha}}|
\underbrace{\widetilde{Z}_{n,\alpha}^{(1)}}_{n\times n_{\alpha}}|\cdots|
\underbrace{\widetilde{Z}_{n,\alpha}^{((n,\alpha)-1)}}_{n\times n_{\alpha}}\right].
\end{eqnarray*}

For $r=0,1,\ldots,n-1$ and $s=0,1,\ldots,n_{\alpha}-1$, we observe
that
\[
(\widetilde{Z}_{n,\alpha}^{(j)})_{r,s}=
(Z_{n,\alpha})_{r,jn_{\alpha}+s},
\]
and
\begin{eqnarray*}
  (Z_{n,\alpha})_{r,jn_{\alpha}+s}&=&\delta_{r-\alpha(jn_{\alpha}+s)}\\
  &=&\delta_{r-j\alpha n_{\alpha}-\alpha s}\\
  &{\underset{\rm (a)}=}&\delta_{r-\alpha s}\\
  &=&(\widetilde{Z}_{n,\alpha}^{(0)})_{r,s}=(\widetilde{Z}_{n,\alpha})_{r,s},
\end{eqnarray*}
where $n_{\alpha}=\frac{n}{(n,\alpha)}$ and (a) is a consequence
of the fact that $\frac{\alpha}{(n,\alpha)}$ is an integer greater
than zero and so $j\alpha
n_{\alpha}=j\frac{\alpha}{(n,\alpha)}n\equiv 0 $ (mod $n$). Thus we
conclude that
$\widetilde{Z}_{n,\alpha}^{(j)}=\widetilde{Z}_{n,\alpha}^{(0)}=\widetilde{Z}_{n,\alpha}$
for $j=0,\ldots,(n,\alpha)-1$.
\end{proof}

%\begin{proof} (of relation $(\ref{3}).$) For $r=0,1,\ldots,n-1$ and $s=0,1,\ldots,n_{\alpha}-1,$ one has
%\begin{eqnarray*}
%    (\widetilde{C}_{n,\alpha})_{r,s}&=&(C_{n,\alpha})_{r,s}=a_{(r-\alpha s)\textrm{ mod}\,n},\\
%    (C_{n})_{r,s}=a_{(r-s)\textrm{ mod}\,n}\\
%    (\widetilde{Z}_{n,\alpha})_{r,s}&=&\delta_{r-\alpha s},
%\end{eqnarray*}
%and
%\begin{eqnarray*}
%  (C_{n}\widetilde{Z}_{n,\alpha})_{r,s} &=& \overset{n-1}
%  {\underset{l=0}\sum}(C_{n})_{r,l}(\widetilde{Z}_{n,\alpha})_{l,s} \\
%    &=& \overset{n-1}{\underset{l=0}\sum}\delta_{l-\alpha s}(C_{n})_{r,l}  \\
%    &=& \overset{n-1}{\underset{l=0}\sum}\delta_{l-\alpha s}a_{(r-l)\textrm{ mod}\,n}  \\
%    &=_{(a)}& a_{(r-\alpha s)\textrm{ mod}\,n} \\
%    &=&  (\widetilde{C}_{n,\alpha})_{r,s},
%\end{eqnarray*}
%where $(a)$ follows from the fact that there exists a unique
%$l\in\{0,1,\ldots,n-1\}$ such that $l-\alpha s\equiv(0)\textrm{ mod
%$n$}$, that is, $l\equiv (\alpha s)\textrm{ mod $n$}$, so,
%$r-l\equiv r-(\alpha s)\textrm{ mod $n$}$ and $(r-l)\textrm{ mod
%$n$}=(r-(\alpha s)\textrm{ mod $n$})\textrm{ mod $n$}=(r-\alpha
%s)\textrm{ mod $n$}$.
%\end{proof}

Another useful fact is represented by the following equation
\begin{eqnarray}\label{3bis}
  \widetilde{Z}_{n,\alpha}=\widetilde{Z}_{n,(n,\alpha)}Z_{n_{\alpha},\check{\alpha}},
\end{eqnarray}
where $Z_{n_{\alpha},\check{\alpha}}$ is the matrix defined in
$(\ref{i})$ of dimension $n_{\alpha}\times n_{\alpha}$. Therefore
\begin{eqnarray}\label{pna}
Z_{n_{\alpha},\check{\alpha}}=\left[\widehat{\delta}_{r-\check{\alpha}
s}
\right]_{r,s=0}^{n_{\alpha}-1},\qquad\widehat{\delta}_{k}=\left\{\begin{array}{cl}
1 & \textrm{if $k\equiv 0 \textrm{ (mod $n_{\alpha}$),}$}\\0 &
\textrm{otherwise.}\end{array} \right.
\end{eqnarray}

Relation $(\ref{3bis})$ will be used later.
\begin{proof} (of relation $(\ref{3bis}).$) For $r=0,1,\ldots,n-1$ and $s=0,1,\ldots,n_{\alpha}-1,$
we find
\begin{eqnarray*}
    (\widetilde{Z}_{n,\alpha})_{r,s}&=&\delta_{r-\alpha s}\\
    &=&\delta_{(r-\alpha s)\textrm{ mod $n$}},
\end{eqnarray*}
and
\begin{eqnarray*}
  (\widetilde{Z}_{n,(n,\alpha)}Z_{n_{\alpha},\check{\alpha}})_{r,s} &=&
        \sum_{l=0}^{n_{\alpha}-1}(\widetilde{Z}_{n,(n,\alpha)})_{r,l}(Z_{n_{\alpha},\check{\alpha}})_{l,s} \\
    &=& \sum_{l=0}^{n_{\alpha}-1}\delta_{r-(n,\alpha) l}\widehat{\delta}_{l-\check{\alpha} s}  \\
    &{\underset{\rm (a)}=}& \delta_{r-(n,\alpha)\cdot(\check{\alpha} s)\textrm{ mod $n_{\alpha}$}}\\
    &=& \delta_{r-(n,\alpha)\cdot\left(\frac{\alpha}{(n,\alpha)} s\right)\textrm{ mod $n_{\alpha}$}} \\
    &{\underset{\rm (b)}=}& \delta_{r-(\alpha s)\textrm{ mod $n$}}\\
    &=& \delta_{(r-(\alpha s)\textrm{ mod $n$})\textrm{ mod $n$}}\\
    &=& \delta_{(r-\alpha s)\textrm{ mod $n$}},
\end{eqnarray*}
where
\begin{itemize}
  \item[(a)] holds true since there exists a unique $l\in\{0,1,\ldots,n_{\alpha}-1\}$ such that
$ l-\check{\alpha} s\equiv0\textrm{ (mod $n_{\alpha}$)}$, that is,
$l\equiv \check{\alpha} s\textrm{ (mod $n_{\alpha}$)}$ and hence
$\delta_{r-(n,\alpha)l}=\delta_{r-(n,\alpha)\cdot(\check{\alpha}
s)\textrm{ mod $n_{\alpha}$}}$;
\item[(b)] is due to the following property: if we have three integer numbers $\rho,\,\theta,$ and $\gamma$, then
\begin{eqnarray*}
  \rho(\theta\textrm{ mod $\gamma)$}=(\rho\theta)\textrm{ mod $\rho\gamma$}.
\end{eqnarray*}
\end{itemize}
\end{proof}

\begin{lemma}\label{amag} If $\alpha\geq n$ then $Z_{n,\alpha}=Z_{n,{\alpha}^\circ}$ where ${\alpha}^\circ$ is the unique integer which satisfies $\alpha=tn+{\alpha}^\circ$ with $0\leq{\alpha}^\circ<n$ and $t\in\mathbb{N}$; $Z_{n,\alpha}$ is defined in $(\ref{i})$.
\end{lemma}

\begin{oss} One can define ${\alpha}^\circ$ by: ${\alpha}^\circ:=\alpha\,\textrm{mod $n$}$.
\end{oss}

\begin{proof} From $(\ref{i})$ we know that
\begin{eqnarray*}
Z_{n,\alpha}=\left[\delta_{r-\alpha c}\right]_{r,c=0}^{n-1},\qquad\delta_{k}=\left\{\begin{array}{cl} 1 & \textrm{if $k\equiv 0 \textrm{ (mod $n$)}$,}\\0 & \textrm{otherwise.}\end{array} \right.
\end{eqnarray*}

For $r,s=0,1,\ldots,n-1$, one has
\begin{eqnarray*}
  (Z_{n,\alpha})_{r,s}=\delta_{r-\alpha s}=\delta_{r-(tn+{\alpha}^\circ)s}=
  \delta_{r-{\alpha}^\circ s}=(Z_{n,{\alpha}^\circ})_{r,s},
\end{eqnarray*}
since $tns\equiv 0\textrm{ (mod $n$)}$. Whence
$Z_{n,\alpha}=Z_{n,{\alpha}^\circ}$.
%thus
%\begin{eqnarray}\label{hat1}
%  &&(Z_{n,\alpha})_{r,c}=1\quad\Longleftrightarrow\quad r-\alpha c\equiv 0 \textrm{ (mod $n$)},\\
%  \label{hat2}&&(Z_{n,\widehat{\alpha}})_{r,c}=1\quad\Longleftrightarrow\quad r-\widehat{\alpha} c\equiv 0 \textrm{ (mod $n$)}.
%\end{eqnarray}
%
%If $\alpha\geq n$, then we can write $\alpha=tn+(\alpha\,\textrm{mod}\,n)$ for some $t\in\mathbb{N}$, and we get
%
%\begin{eqnarray*}
%  &&r-\alpha c\equiv 0 \textrm{ (mod $n$)}\\
%  &&\qquad\Updownarrow\quad{\scriptstyle \alpha=tn+(\alpha\,\textrm{mod}\,n)} \\
%  &&r-[tn+(\alpha\,\textrm{mod}\,n)c]\equiv 0\textrm{ (mod $n$)}\\
%  &&\qquad\Updownarrow\\
%  &&r-tn-(\alpha\,\textrm{mod}\,n)c\equiv 0 \textrm{ (mod $n$)}\\
%  &&\qquad\Updownarrow\quad {\scriptstyle tn\equiv 0\textrm{ (mod $n$)}}\\
%  &&r-(\alpha\,\textrm{mod}\,n)c\equiv 0 \textrm{ (mod $n$)}\\
%  &&\qquad\Updownarrow\quad{\scriptstyle\widehat{\alpha}=\alpha\,\textrm{ mod $n$}}\\
%  &&r-\widehat{\alpha} c\equiv 0 \textrm{ (mod $n$)},
%\end{eqnarray*}
%
%then the two expressions in $(\ref{hat1})$ and $(\ref{hat2})$ are equivalent.
\end{proof}

The previous lemma tells us that, for $\alpha$-circulant matrices, we can consider only the case where $0\leq\alpha<n$.
In fact, if $\alpha\geq n$, from $(\ref{0})$ we infer that
\begin{eqnarray*}
  C_{n,\alpha}=C_{n}Z_{n,\alpha}=C_{n}Z_{n,{\alpha}^\circ}=C_{n,{\alpha}^\circ}.
\end{eqnarray*}

Finally, it is worth noticing that the use of $(\ref{iV})$ and
$(\ref{0})$ implies that
\begin{eqnarray}\label{Vi}
  C_{n,\alpha} = F_{n}D_{n}F_{n}^{\ast}Z_{n,\alpha}.
\end{eqnarray}

Formula $(\ref{Vi})$ plays an important role for studying the
singular values of the $\alpha$-circulant matrices.

\subsection{A characterization of $Z_{n,\alpha}$ in terms of Fourier matrices}

\begin{lemma}\label{l1}
Let $F_{n}$ be the Fourier matrix of order $n$ defined in (\ref{V})
and let $\widetilde{Z}_{n,\alpha}\in \mathbb{C}^{n\times
n_{\alpha}}$ be the matrix represented in $(\ref{z-n-alfa})$. Then
\begin{eqnarray}\label{4}
  F_{n}\widetilde{Z}_{n,\alpha} =
  \frac{1}{\sqrt{(n,\alpha)}}I_{n,\alpha}F_{n_{\alpha}}Z_{n_{\alpha},\check{\alpha}},
\end{eqnarray}
where $I_{n,\alpha}\in \mathbb{C}^{n\times n_{\alpha}}$ and
\begin{eqnarray*}
  I_{n,\alpha}=\left.\left[\begin{array}{c}
  I_{n_{\alpha}}\\
  \hline
  I_{n_{\alpha}}\\
  \hline
  \vdots\\
  \hline
  I_{n_{\alpha}}
  \end{array}\right]\right\}\textrm{$(n,\alpha)$ times,}
\end{eqnarray*}
with $I_{n_{\alpha}}$ being the identity matrix of size $n_{\alpha}$
and $Z_{n_{\alpha},\check{\alpha}}$ as in $(\ref{pna})$.
\end{lemma}

\begin{oss}
\text{\,\,\,}  $n= n_{\alpha}\cdot(n,\alpha).$
\end{oss}

\begin{proof} (of Lemma $\ref{l1}$.)
Rewrite the Fourier matrix as
\begin{eqnarray*}
  F_{n}=\frac{1}{\sqrt{n}}\left[\begin{array}{c|c|c|c|c}
  f_{0} & f_{1} & f_{2} & \cdots & f_{n-1}
  \end{array}\right],
\end{eqnarray*}
where $f_{k},$ $k=0,1,2,\ldots,n-1,$ is the $k-th$ column of the
Fourier matrix of order $n$:
\begin{eqnarray}\label{6}
  f_{k}=\left[e^{-\frac{2\pi ikj}{n}}\right]_{j=0}^{n-1}=\left[\begin{array}{c}
  e^{-\frac{2\pi ik\cdot0}{n}}\\
  e^{-\frac{2\pi ik\cdot1}{n}}\\
  e^{-\frac{2\pi ik\cdot2}{n}}\\
  \vdots\\
  e^{-\frac{2\pi ik\cdot(n-1)}{n}}\\
  \end{array}\right].
\end{eqnarray}

From $(\ref{3bis})$, we find
\begin{equation}\label{7}
  F_{n}\widetilde{Z}_{n,\alpha}=F_{n}\widetilde{Z}_{n,(n,\alpha)}Z_{n_{\alpha},\check{\alpha}}=
  \frac{1}{\sqrt{n}}\left[\begin{array}{c|c|c|c|c}
  f_{0} & f_{1\cdot(n,\alpha)} & f_{2\cdot(n,\alpha)} & \cdots & f_{(n_{\alpha}-1)\cdot(n,\alpha)}
  \end{array}\right]Z_{n_{\alpha},\check{\alpha}}\in \mathbb{C}^{n\times
  n_{\alpha}}.
\end{equation}

Indeed, for $k=0,1,\ldots,n_{\alpha}-1$, $j=0,1,\ldots,n-1,$ one has
\begin{equation}\label{8}
    \left(F_{n}\widetilde{Z}_{n,(n,\alpha)}\right)_{j,k}=\overset{n-1}
    {\underset{l=0}\sum}(F_{n})_{j,l}(\widetilde{Z}_{n,(n,\alpha)})_{l,k}=
    \overset{n-1}{\underset{l=0}\sum}
    \delta_{l-(n,\alpha)k}e^{-\frac{2\pi ijl}{n}},
\end{equation}
and, since $0\leq (n,\alpha)k\leq n-(n,\alpha),$ there exists a
unique $l_{k}\in \{0,1,2,\ldots,n-1\}$ such that
$l_{k}-(n,\alpha)k\equiv0$ (mod $n$), so $l_{k}=(n,\alpha)k$.
Consequently relation $(\ref{8})$ implies
\begin{equation*}
  \left(F_{n}\widetilde{Z}_{n,(n,\alpha)}\right)_{j,k}=\delta_{l_{k}-(n,\alpha)k}
  e^{-\frac{2\pi ijl_{k}}{n}}=e^{-\frac{2\pi ij
  (n,\alpha)k}{n}}=\left(f_{(n,\alpha)k}\right)_{j},
  \end{equation*}
for all $0\leq j\leq n-1$ and $0\leq k\leq n_{\alpha}-1,$ and hence
\begin{equation*}
    F_{n}\widetilde{Z}_{n,(n,\alpha)}=
  \frac{1}{\sqrt{n}}\left[\begin{array}{c|c|c|c|c}
  f_{0} & f_{1\cdot(n,\alpha)} & f_{2\cdot(n,\alpha)} & \cdots & f_{(n_{\alpha}-1)\cdot(n,\alpha)}
  \end{array}\right].
\end{equation*}

For $k=0,1,2,\ldots,n_{\alpha}-1,$ we deduce
\begin{equation*}
    f_{(n,\alpha)k}=\left[e^{-\frac{2\pi ij(n,\alpha)k}{n}}\right]_{j=0}^{n-1}=\left[e^{-\frac{2\pi ijk}{n_{\alpha}}}\right]_{j=0}^{n-1},
\end{equation*}
and then, taking into account the equalities
$n=(n,\alpha)\frac{n}{(n,\alpha)}=(n,\alpha)n_{\alpha},$ we can
write
\begin{eqnarray}\label{9}
  f_{(n,\alpha)k} =\left[\begin{array}{l}
   \left[e^{-\frac{2\pi ikj}{n_{\alpha}}}\right]_{j=0}^{n_{\alpha}-1} \\
   \left[e^{-\frac{2\pi ikj}{n_{\alpha}}}\right]_{j=n_{\alpha}}^{2n_{\alpha}-1} \\
               \qquad\;\vdots            \\
   \left[e^{-\frac{2\pi ikj}{n_{\alpha}}}\right]_{j=((n,\alpha)-1)n_{\alpha}}^{(n,\alpha)
   n_{\alpha}-1}
   \end{array}\right],
\end{eqnarray}
where
\begin{eqnarray}\label{10}
  \left[e^{-\frac{2\pi ikj}{n_{\alpha}}}\right]_{j=0}^{n_{\alpha}-1} =
   \begin{bmatrix}
  e^{-\frac{2\pi ik\cdot0}{n_{\alpha}}} \\
  e^{-\frac{2\pi ik\cdot1}{n_{\alpha}}} \\
  e^{-\frac{2\pi ik\cdot2}{n_{\alpha}}} \\
          \vdots                  \\
  e^{-\frac{2\pi ik\cdot(n_{\alpha}-1)}{n_{\alpha}}}
  \end{bmatrix}.
\end{eqnarray}

According to formula $(\ref{6}),$ one observes that the vector in
$(\ref{10})$ is the $k-th$ column of the Fourier matrix
$F_{n_{\alpha}}$. Furthermore, for $l=0,1,2,\ldots,(n,\alpha)-1$, we
find
\begin{equation}\label{11}
  \left[e^{-\frac{2\pi ikj}{n_{\alpha}}}\right]_{j=ln_{\alpha}}^{(l+1)n_{\alpha}-1}
      =
  \begin{bmatrix}
   e^{-\frac{2\pi ikln_{\alpha}}{n_{\alpha}}} \\
   e^{-\frac{2\pi ik(ln_{\alpha}+1)}{n_{\alpha}}} \\
   e^{-\frac{2\pi ik(ln_{\alpha}+2)}{n_{\alpha}}} \\
                  \vdots                                \\
   e^{-\frac{2\pi ik(ln_{\alpha}+n_{\alpha}-1)}{n_{\alpha}}}
 \end{bmatrix}
                                      =e^{-2\pi ikl}
 \begin{bmatrix}
 e^{-\frac{2\pi ik\cdot0}{n_{\alpha}}} \\
 e^{-\frac{2\pi ik\cdot1}{n_{\alpha}}} \\
 e^{-\frac{2\pi ik\cdot2}{n_{\alpha}}} \\
                \vdots                   \\
 e^{-\frac{2\pi ik\cdot(n_{\alpha}-1)}{n_{\alpha}}}
 \end{bmatrix}
   =\left[e^{-\frac{2\pi ikj}{n_{\alpha}}}\right]_{j=0}^{n_{\alpha}-1}.
\end{equation}

Using $(\ref{11})$, the expression of the vector in $(\ref{9})$
becomes
 \begin{eqnarray}\label{12}
 f_{(n,\alpha)k}=\left.
 \begin{bmatrix}
  \left[e^{-\frac{2\pi ikj}{n_{\alpha}}}\right]_{j=0}^{n_{\alpha}-1} \\
  \left[e^{-\frac{2\pi ikj}{n_{\alpha}}}\right]_{j=0}^{n_{\alpha}-1} \\
              \!\!\!\!\!\vdots             \\
  \left[e^{-\frac{2\pi ikj}{n_{\alpha}}}\right]_{j=0}^{n_{\alpha}-1}
 \end{bmatrix}
 \right\}\textrm{ $(n,\alpha)$ times.}
\end{eqnarray}

Setting $\widetilde{f}_{r}=\left[e^{-\frac{2\pi irj}
 {n_{\alpha}}}\right]_{j=0}^{n_{\alpha}-1},$
for $0\leq r\leq n_{\alpha}-1$, the Fourier matrix $F_{n_{\alpha}}$
of size $n_{\alpha}$ takes the form
\begin{eqnarray}\label{13}
  F_{n_{\alpha}}=\frac{1}{\sqrt{n_{\alpha}}}\left[\begin{array}{c|c|c|c|c}
  \widetilde{f}_{0} & \widetilde{f}_{1} & \widetilde{f}_{2} & \cdots & \widetilde{f}_{n_{\alpha}-1}
  \end{array}\right].
\end{eqnarray}

From formula $(\ref{10}),$ the relation $(\ref{12})$ can be
expressed as
\begin{eqnarray*}
  f_{(n,\alpha)k}=\left.\left[\begin{array}{c}
  \widetilde{f}_{k}\\
  \widetilde{f}_{k}\\
  \vdots\\
  \widetilde{f}_{k}
  \end{array}\right]\right\}\textrm{ $(n,\alpha)$ times,} \qquad k=0,\ldots,n_{\alpha}-1,
\end{eqnarray*}
and, as a consequence, formula $(\ref{7})$ can be rewritten as
\begin{eqnarray*}
  F_{n}\widetilde{Z}_{n,\alpha}=F_{n}\widetilde{Z}_{n,(n,\alpha)}Z_{n_{\alpha},\check{\alpha}} &=&
  \frac{1}{\sqrt{n}}\left[\begin{array}{c|c|c|c|c}
  \widetilde{f}_{0} & \widetilde{f}_{1} & \widetilde{f}_{2} & \cdots & \widetilde{f}_{n_{\alpha}-1}\\
  \widetilde{f}_{0} & \widetilde{f}_{1} & \widetilde{f}_{2} & \cdots & \widetilde{f}_{n_{\alpha}-1}\\
  \vdots & \vdots & \vdots & \vdots & \vdots \\
  \widetilde{f}_{0} & \widetilde{f}_{1} & \widetilde{f}_{2} & \cdots & \widetilde{f}_{n_{\alpha}-1}
  \end{array}\right]Z_{n_{\alpha},\check{\alpha}}\\
  &=&\frac{1}{\sqrt{(n,\alpha)n_{\alpha}}}\left[\begin{array}{c}
  \sqrt{n_{\alpha}}F_{n_{\alpha}}\\
  \hline
  \sqrt{n_{\alpha}}F_{n_{\alpha}}\\
  \hline
  \vdots  \\
  \hline
  \sqrt{n_{\alpha}}F_{n_{\alpha}}
  \end{array}\right]Z_{n_{\alpha},\check{\alpha}}\\
  &=&\frac{1}{\sqrt{(n,\alpha)}}\left[\begin{array}{c}
  F_{n_{\alpha}}\\
  \hline
  F_{n_{\alpha}}\\
  \hline
  \vdots  \\
  \hline
  F_{n_{\alpha}}
  \end{array}\right]Z_{n_{\alpha},\check{\alpha}}\\
  &=&\frac{1}{\sqrt{(n,\alpha)}}\left[\begin{array}{c}
  I_{n_{\alpha}}\\
  \hline
  I_{n_{\alpha}}\\
  \hline
  \vdots  \\
  \hline
  I_{n_{\alpha}}
  \end{array}\right]F_{n_{\alpha}}Z_{n_{\alpha},\check{\alpha}}\\
  &=& \frac{1}{\sqrt{(n,\alpha)}}I_{n,\alpha}F_{n_{\alpha}}Z_{n_{\alpha},\check{\alpha}}.
\end{eqnarray*}
\end{proof}

In the subsequent subsection, we will exploit Lemma $\ref{l1}$ in
order to characterize the singular values of the $\alpha$-circulant
matrices $C_{n,\alpha}$. Here we conclude the subsection with the
following simple observations.
\begin{oss} In Lemma $\ref{l1}$, if $(n,\alpha)=\alpha$, we
have  $n_{\alpha}=\frac{n}{(n,\alpha)}=\frac{n}{\alpha}$ and
$\check{\alpha}=\frac{\alpha}{(n,\alpha)}=1$; so the matrix
$Z_{n_{\alpha},\check{\alpha}}=Z_{n_{\alpha},1}$, appearing in
$(\ref{4})$, is the identity matrix of dimension
$\frac{n}{\alpha}\times\frac{n}{\alpha}$. The relation $(\ref{4})$
becomes
\begin{eqnarray*}
  F_{n}\widetilde{Z}_{n,\alpha}=\frac{1}{\sqrt{\alpha}}I_{n,\alpha}F_{n_{\alpha}}.
\end{eqnarray*}

The latter equation with $\alpha=2$ and even $n$ appear (and is
crucial) in the multigrid literature; see {\rm \cite{mcirco}},
equation (3.2), page 59 and, in slightly different form for the sine
algebra of type I, see {\rm \cite{FS1}}, Section 2.1.
\end{oss}

\begin{oss} If $(n,\alpha)=1$, Lemma $\ref{l1}$ is trivial, because
$n_{\alpha}=\frac{n}{(n,\alpha)}=n$,
$\check{\alpha}=\frac{\alpha}{(n,\alpha)}=\alpha$, and so
$\widetilde{Z}_{n,\alpha}=Z_{n,\alpha}$. The relation $(\ref{4})$
becomes
 \begin{eqnarray*}
    F_{n}\widetilde{Z}_{n,\alpha}=F_{n}Z_{n,\alpha}&=&I_{n,\alpha}F_{n_{\alpha}}Z_{n_{\alpha},\check{\alpha}}\\
    &=&F_{n}Z_{n,\alpha},
  \end{eqnarray*}
since the matrix $I_{n,\alpha}$ reduces by its definition to the
identity matrix of order $n$.
\end{oss}

\begin{oss}
Lemma $\ref{l1}$ is true also if, instead of $F_{n}$ and
$F_{n_{\alpha}}$, we put $F_{n}^{*}$ and $F_{n_{\alpha}}^{*}$,
respectively, because $F_{n}^{*}=\overline{F_{n}}$. In fact there is
no transposition, but only conjugation.
\end{oss}

\subsection{Characterization of the singular values of the $\alpha$-circulant matrices}

Now we link the singular values of $\alpha$-circulant matrices with
the eigenvalues of its circulant counterpart $C_n$. This is
nontrivial given the multiplicative relation $C_{n,\alpha}= C_n
Z_{n,\alpha}$.

Having in mind the definition of the diagonal matrix $D_n$ given in
(\ref{V}), we start by setting
\begin{eqnarray}\label{DJ}
  \notag&&D_{n}^{\ast}D_{n}={\rm diag}(|D_{n}|^{2}_{s,s};\,\,s=0,1,\ldots,n-1)=
    {\rm diag}(d_{s};\,\,s=0,1,\ldots,n-1)=
    \overset{(n,\alpha)}{\underset{l=1}\oplus}\Delta_{l},\\
    &&J_{(n,\alpha)}\otimes I_{n_{\alpha}} =\underset{(n,\alpha)\text{\,\,}times}{
    \underbrace{\left[I_{n,\alpha}|I_{n,\alpha}|\cdots|I_{n,\alpha}\right]}}=
\left.\left[\begin{array}{c|c|c|c}
I_{n_{\alpha}}  & I_{n_{\alpha}} & \cdots & I_{n_{\alpha}}\\
\hline
I_{n_{\alpha}}  & I_{n_{\alpha}} & \cdots & I_{n_{\alpha}}\\
\hline
\vdots &  \vdots & \vdots & \vdots\\
\hline
I_{n_{\alpha}}  & I_{n_{\alpha}} & \cdots & I_{n_{\alpha}}
\end{array}\right]\right\}\textrm{ $(n,\alpha)$ times,}
\end{eqnarray}
where
\begin{eqnarray}
   \label{ds}&&d_{s}=|D_{n}|_{s,s}^{2}=(D_{n})_{s,s}\cdot\overline{(D_{n})_{s,s}},
   \quad\textrm{$D_{n}$ defined in $(\ref{V})$, $s=0,1,\ldots,n-1$,}\\
   \notag&&\Delta_{l}=\left[\begin{array}{cccc}
   d_{(l-1)n_{\alpha}} & & & \\
   & d_{(l-1)n_{\alpha}+1} & & \\
   & & \ddots & \\
   & & & d_{(l-1)n_{\alpha}+n_{\alpha}-1}
   \end{array}\right]\in \mathbb{C}^{n_{\alpha}\times
   n_{\alpha}};\text{\,\,}l=1,2,\ldots,(n,\alpha),\\
   \label{19bis}&& J_{(n,\alpha)} =\left.\left[\begin{array}{cccc}
   1 & 1 & \cdots & 1\\
   1 & 1 & \cdots & 1\\
   \vdots & \vdots & \vdots & \vdots\\
   1 & 1 & \cdots & 1
\end{array}\right]
   \right\}\textrm{ $(n,\alpha)$ times.}
\end{eqnarray}

We now exploit relation $(\ref{18})$ and Lemma $\ref{l1}$, and we
obtain that
  \begin{eqnarray}\label{fnz}
  \notag F_{n}Z_{n,\alpha} &=& F_{n}\left[\widetilde{Z}_{n,\alpha}|\widetilde{Z}_{n,\alpha}
  |\cdots|\widetilde{Z}_{n,\alpha}\right] \\
    \notag&=& \left[F_{n}\widetilde{Z}_{n,\alpha}|F_{n}\widetilde{Z}_{n,\alpha}|\cdots
    |F_{n}\widetilde{Z}_{n,\alpha}\right] \\
    \notag&=& \frac{1}{\sqrt{(n,\alpha)}}\left[I_{n,\alpha}F_{n_{\alpha}}Z_{n_{\alpha},\check{\alpha}}
    |I_{n,\alpha}F_{n_{\alpha}}Z_{n_{\alpha},\check{\alpha}}|\cdots|I_{n,\alpha}F_{n_{\alpha}}Z_{n_{\alpha},\check{\alpha}}\right] \\
    \notag&=& \frac{1}{\sqrt{(n,\alpha)}}\left[I_{n,\alpha}|I_{n,\alpha}|\cdots|I_{n,\alpha}\right]
    \left.\left[\begin{array}{cccc}
    F_{n_{\alpha}}Z_{n_{\alpha},\check{\alpha}} &   &   &      \\
    &\!\!\!\!\! F_{n_{\alpha}}Z_{n_{\alpha},\check{\alpha}} &   &     \\
    &   &    \!\!\!\!\!\ddots    &   \\
    &   &   &   \!\!\!\!\!F_{n_{\alpha}}Z_{n_{\alpha},\check{\alpha}} \\
    \end{array}\right]\right\}\textrm{ $(n,\alpha)$ times}\\
    &=&\frac{1}{\sqrt{(n,\alpha)}}\left[I_{n,\alpha}|I_{n,\alpha}|\cdots
    |I_{n,\alpha}\right]\left(I_{(n,\alpha)}\otimes F_{n_{\alpha}}Z_{n_{\alpha},\check{\alpha}}\right),
  \end{eqnarray}
where $I_{(n,\alpha)}$ is the identity matrix of order $(n,\alpha).$
  Furthermore,
 \begin{eqnarray}\label{CastC}
  \notag C_{n,\alpha}^{\ast}C_{n,\alpha} &=&(F_{n}D_{n}F_{n}^{\ast}Z_{n,\alpha})^{\ast}(F_{n}D_{n}F_{n}^{\ast}Z_{n,\alpha}) \\
                                  &=&
\notag Z_{n,\alpha}^{\ast}F_{n}D_{n}^{\ast}F_{n}^{\ast}F_{n}D_{n}F_{n}^{\ast}Z_{n,\alpha} \\
                                  &=&
\notag Z_{n,\alpha}^{\ast}F_{n}D_{n}^{\ast}D_{n}F_{n}^{\ast}Z_{n,\alpha} \\
                                  &=&
(F_{n}^{\ast}Z_{n,\alpha})^{\ast}D_{n}^{\ast}D_{n}F_{n}^{\ast}Z_{n,\alpha}.
\end{eqnarray}

From $(\ref{fnz})$ and $(\ref{DJ})$, we plainly infer the following
relations
\begin{eqnarray*}
  (F_{n}^{\ast}Z_{n,\alpha})^{\ast}&=&\left(\frac{1}{\sqrt{(n,\alpha)}}\left[I_{n,\alpha}|I_{n,\alpha}|\cdots
|I_{n,\alpha}\right]\left(I_{(n,\alpha)}\otimes F_{n_{\alpha}}^{*}Z_{n_{\alpha},\check{\alpha}}\right)\right)^{*}\\
&=&\frac{1}{\sqrt{(n,\alpha)}}\left(I_{(n,\alpha)}\otimes F_{n_{\alpha}}^{*}Z_{n_{\alpha},\check{\alpha}}\right)^{*}\left(J_{(n,\alpha)}\otimes I_{n_{\alpha}}\right)\\
&=&\frac{1}{\sqrt{(n,\alpha)}}\left(I_{(n,\alpha)}\otimes Z_{n_{\alpha},\check{\alpha}}^{*}F_{n_{\alpha}}\right)\left(J_{(n,\alpha)}\otimes I_{n_{\alpha}}\right),\\
F_{n}^{\ast}Z_{n,\alpha}&=&\frac{1}{\sqrt{(n,\alpha)}}\left[I_{n,\alpha}|I_{n,\alpha}|\cdots
|I_{n,\alpha}\right]\left(I_{(n,\alpha)}\otimes F_{n_{\alpha}}^{*}Z_{n_{\alpha},\check{\alpha}}\right)\\
&=&\frac{1}{\sqrt{(n,\alpha)}}\left(J_{(n,\alpha)}\otimes
I_{n_{\alpha}}\right)\left(I_{(n,\alpha)}\otimes
F_{n_{\alpha}}^{*}Z_{n_{\alpha},\check{\alpha}}\right).
\end{eqnarray*}

Hence
\begin{eqnarray*}
C_{n,\alpha}^{\ast}C_{n,\alpha}=\!\left(I_{(n,\alpha)}\otimes
Z_{n_{\alpha},\check{\alpha}}^{*}F_{n_{\alpha}}\right)\!\left(J_{(n,\alpha)}\otimes
I_{n_{\alpha}}\right)\!\frac{1}{(n,\alpha)}D_{n}^{\ast}D_{n}\left(J_{(n,\alpha)}\otimes
I_{n_{\alpha}}\right)\!\left(I_{(n,\alpha)} \otimes
F_{n_{\alpha}}^{*}Z_{n_{\alpha},\check{\alpha}}\right).
\end{eqnarray*}

Now using the properties of the tensorial product
\begin{eqnarray*}
  (I_{(n,\alpha)}\otimes Z_{n_{\alpha},\check{\alpha}}^{*}F_{n_{\alpha}})\!\!\!\!\!\!\!\!&&\!\!\!\!\!\!\!\!(I_{(n,\alpha)}\otimes F_{n_{\alpha}}^{*}Z_{n_{\alpha},\check{\alpha}})\\
  &=&I_{(n,\alpha)}I_{(n,\alpha)}\otimes Z_{n_{\alpha},\check{\alpha}}^{*}F_{n_{\alpha}}F_{n_{\alpha}}^{*}Z_{n_{\alpha},\check{\alpha}}\\
  &=&I_{(n,\alpha)}I_{(n,\alpha)}\otimes Z_{n_{\alpha},\check{\alpha}}^{*}Z_{n_{\alpha},\check{\alpha}}\\
  &=&I_{(n,\alpha)}I_{(n,\alpha)}\otimes I_{n_{\alpha}}=I_{n},
\end{eqnarray*}
and from a similarity argument, one deduces  that the eigenvalues of
$C_{n,\alpha}^{\ast}C_{n,\alpha}$ are the eigenvalues of the matrix
\begin{eqnarray*}
&&\!\!\!\!\!\!\!\!\!\!\!\!\!\!\!\!\left(J_{(n,\alpha)}\otimes I_{n_{\alpha}}\right)\frac{1}{(n,\alpha)}D_{n}^{\ast}D_{n}\left(J_{(n,\alpha)}\otimes I_{n_{\alpha}}\right)\\
&=&\frac{1}{(n,\alpha)}\left[\begin{array}{c|c|c|c}
I_{n_{\alpha}}  & I_{n_{\alpha}} & \cdots & I_{n_{\alpha}}\\
\hline
I_{n_{\alpha}}  & I_{n_{\alpha}} & \cdots & I_{n_{\alpha}}\\
\hline
\vdots &  \vdots & \vdots & \vdots\\
\hline
I_{n_{\alpha}}  & I_{n_{\alpha}} & \cdots & I_{n_{\alpha}}
\end{array}\right]
\left[\begin{array}{cccc}
  \Delta_{1} &  &  &  \\
  & \Delta_{2} &  &   \\
  &          & \ddots &  \\
  &            &  &  \Delta_{(n,\alpha)}
\end{array}\right]
\left[\begin{array}{c|c|c|c}
I_{n_{\alpha}}  & I_{n_{\alpha}} & \cdots & I_{n_{\alpha}}\\
\hline
I_{n_{\alpha}}  & I_{n_{\alpha}} & \cdots & I_{n_{\alpha}}\\
\hline
\vdots &  \vdots & \vdots & \vdots\\
\hline
I_{n_{\alpha}}  & I_{n_{\alpha}} & \cdots & I_{n_{\alpha}}
\end{array}\right]\\
&=&\frac{1}{(n,\alpha)}\left[\begin{array}{c|c|c|c}
I_{n_{\alpha}}  & I_{n_{\alpha}} & \cdots & I_{n_{\alpha}}\\
\hline
I_{n_{\alpha}}  & I_{n_{\alpha}} & \cdots & I_{n_{\alpha}}\\
\hline
\vdots &  \vdots  &\vdots & \vdots\\
\hline
I_{n_{\alpha}}  & I_{n_{\alpha}} &\cdots & I_{n_{\alpha}}
\end{array}\right]
\left[\begin{array}{c|c|c|c}
  \Delta_{1} & \Delta_{1} & \cdots & \Delta_{1}  \\
  \hline
  \Delta_{2} & \Delta_{2} & \cdots & \Delta_{2}  \\
  \hline
  \vdots & \vdots & \vdots & \vdots \\
  \hline
  \Delta_{(n,\alpha)} & \Delta_{(n,\alpha)} & \cdots & \Delta_{(n,\alpha)}
\end{array}\right]\\
&=&\frac{1}{(n,\alpha)}\left[\begin{array}{c|c|c|c}
  \overset{(n,\alpha)}{\underset{l=1}\sum}\Delta_{l} & \overset{(n,\alpha)}{\underset{l=1}\sum}\Delta_{l} & \cdots & \overset{(n,\alpha)}{\underset{l=1}\sum}\Delta_{l}  \\
  \hline
  \overset{(n,\alpha)}{\underset{l=1}\sum}\Delta_{l} & \overset{(n,\alpha)}{\underset{l=1}\sum}\Delta_{l} & \cdots & \overset{(n,\alpha)}{\underset{l=1}\sum}\Delta_{l}  \\
  \hline
  \vdots & \vdots & \vdots & \vdots \\
  \hline
  \overset{(n,\alpha)}{\underset{l=1}\sum}\Delta_{l} & \overset{(n,\alpha)}{\underset{l=1}\sum}\Delta_{l} & \cdots & \overset{(n,\alpha)}{\underset{l=1}\sum}\Delta_{l}
\end{array}\right]\\
&=&\frac{1}{(n,\alpha)}\underbrace{\left[\begin{array}{cccc}
   1 & 1 & \cdots & 1\\
   1 & 1 & \cdots & 1\\
   \vdots & \vdots & \vdots & \vdots\\
   1 & 1 & \cdots & 1
\end{array}\right]}_{\textrm{$(n,\alpha)$ times}}\otimes\left(\sum_{l=1}^{(n,\alpha)}\Delta_{l}\right).
\end{eqnarray*}

Therefore, from $(\ref{19bis})$, we infer that
\begin{equation}\label{eigg}
    {\rm Eig}(C_{n,\alpha}^{\ast}C_{n,\alpha})=\frac{1}{(n,\alpha)}{\rm Eig}\left(J_{(n,\alpha)}
    \otimes \overset{(n,\alpha)}{\underset{l=1}\sum}\Delta_{l}\right),
\end{equation}
where
\begin{eqnarray}\label{19}
  \frac{1}{(n,\alpha)}{\rm Eig}(J_{(n,\alpha)}) = \{0,1\}.
\end{eqnarray}

Here we must observe that $\frac{1}{(n,\alpha)}J_{(n,\alpha)}$ is a
matrix of rank 1, so it has all eigenvalues equal to zero except one
eigenvalue equal to 1. In fact note that the trace of a matrix is,
by definition, the sum of its eigenvalues: in our case the trace is
$(n,\alpha)\cdot\frac{1}{(n,\alpha)}=1$ and hence the only nonzero
eigenvalue is necessarily equal to 1. Moreover
\begin{eqnarray*}
    \overset{(n,\alpha)}{\underset{l=1}\sum}\Delta_{l}&=&
    \overset{(n,\alpha)}{\underset{l=1}\sum}{\rm diag}(d_{(l-1)n_{\alpha}+j};\text{\,\,}j=0,1,\ldots,n_{\alpha}-1)\\
    &=&{\rm diag}\left(\overset{(n,\alpha)}{\underset{l=1}\sum}d_{(l-1)n_{\alpha}+j};\text{\,\,}
    j=0,1,\ldots,n_{\alpha}-1\right).
    \end{eqnarray*}

Consequently, since
$\overset{(n,\alpha)}{\underset{l=1}\sum}\Delta_{l}$ is a diagonal
matrix, we have
    \begin{eqnarray}\label{20}
        {\rm Eig}\left(\overset{(n,\alpha)}{\underset{l=1}\sum}\Delta_{l}\right)&=&
\left\{\overset{(n,\alpha)}{\underset{l=1}\sum}d_{(l-1)n_{\alpha}+j};\text{\,\,}
    j=0,1,\ldots,n_{\alpha}-1\right\},
    \end{eqnarray}
where $d_{k}$ are defined in $(\ref{ds})$.

Finally, by exploiting basic properties of the tensor product, we
know that the eigenvalues of a tensor product of two square matrices
$A\otimes B$ are given by all possible products of  eigenvalues of
$A$ of order $p$ and of eigenvalues of $B$ of order $q$, that is
$\lambda(A\otimes B)=\lambda_{j}(A)\lambda_{k}(B)$ for
$j=1,\ldots,p$ and $k=1,\ldots,q$. Therefore, by taking into
consideration $(\ref{eigg})$, $(\ref{19})$, and $(\ref{20})$, we
find
\begin{eqnarray}\label{21}
  \lambda_{j}(C_{n,\alpha}^{\ast}C_{n,\alpha})&=&\overset{(n,\alpha)}
  {\underset{l=1}\sum}d_{(l-1)n_{\alpha}+j}, \quad j=0,1,\ldots,n_{\alpha}-1,\\
  \label{21bis}\lambda_{j}(C_{n,\alpha}^{\ast}C_{n,\alpha})&=& 0,\qquad j=n_{\alpha},\ldots,n-1.
\end{eqnarray}

From $(\ref{21})$, $(\ref{21bis})$ and $(\ref{1}),$ one obtains that
the singular values of an $\alpha$-circulant matrix $C_{n,\alpha}$
are given by
\begin{eqnarray}\label{22}
  \sigma_{j}(C_{n,\alpha})&=&\sqrt{\overset{(n,\alpha)}
  {\underset{l=1}\sum}d_{(l-1)n_{\alpha}+j}}, \quad j=0,1,\ldots,n_{\alpha}-1,\\
  \notag\sigma_{j}(C_{n,\alpha})&=& 0,\qquad j=n_{\alpha},\ldots,n-1,
\end{eqnarray}
where the values $d_{k}$, $k=0,\ldots,n-1$, are defined in
$(\ref{ds})$.

\subsection{Special cases and observations}\label{special-case}

In this subsection we consider some special cases and we furnish a
further link between the eigenvalues of circulant matrices and the
singular values of $\alpha$-circulants. In the case where
$(n,\alpha)=1,$ we have $n_{\alpha}=\frac{n}{(n,\alpha)}=n$. Hence
the formula $(\ref{22})$ becomes
\begin{eqnarray*}
  \sigma_{j}(C_{n,\alpha})=\sqrt{d_{j}}, \quad j=0,1,\ldots,n-1.
\end{eqnarray*}

In other words the singular values of $C_{n,\alpha}$ coincide with
those of $C_n$ (this is expected since $Z_{n,\alpha}$ is a
permutation matrix) and in particular with the moduli of the
eigenvalues of $C_n$.

Concerning the eigenvalues of circulant matrices it should be
observed that formula (\ref{V}) can be interpreted in function terms
as the evaluation of a polynomial at the grid points given by the
$n$-th roots of the unity. This is a standard observation because
the Fourier matrix is a special instance of the classical
Vandermonde matrices when the knots are exactly all the $n$-th roots
of the unity.

Therefore, defining the polynomial $p(t)=\sum_{k=0}^{n-1} a_k e^{ikt}$,
it is trivial to observe that the eigenvalues of $C_n= F_n D_n
F_n^*$ are given by
\[
\lambda_j(C_n)=p\left(\frac{2\pi j}{n}\right), \ \ \ j=0,\ldots,n-1.
\]

The question that naturally arises is how to connect the expression
in (\ref{22}) of the nontrivial singular values of $C_{n,\alpha}$
with the polynomial $p$. The answer is somehow intriguing and can be
resumed in the following formula which could be of interest in the
multigrid community (see Section \ref{sec:multigrid})
\begin{eqnarray}\label{sv-nontrivial-symbol-p}
  \sigma_{j}(C_{n,\alpha}) & = & \sqrt{\overset{(n,\alpha)-1}
  {\underset{l=0}\sum} |p|^2\left(\frac{x_j+2\pi l}{(n,\alpha)} \right)},
  \quad x_j=\frac{2\pi j}{n_{\alpha}}, \ \ j=0,1,\ldots,n_{\alpha}-1.
\end{eqnarray}

In addition if $\alpha$ is fixed and a sequence of integers $n$ is
chosen so that $(n,\alpha)>1$ for $n$ large enough, then
$\{C_{n,\alpha}\}\sim_\sigma (0,G)$ for a proper set $G$. If the
sequence of $n$ is chosen so that $n$ and $\alpha$ are coprime for
all $n$ large enough, then the existence of the distribution is
related to the smoothness properties of a function $f$ such that
$\{a_k\}$ can be interpreted as the sequence of its Fourier
coefficients (see e.g. \cite{appr-mult}). From the above reasoning
it is clear that, if $n$ is allowed to be vary among all the
positive integer numbers, then $\{C_{n,\alpha}\}$ does not possess a
joint singular value distribution.

\section{Singular values of $\alpha$-Toeplitz matrices}\label{sec:toep}

For $p=q=d=1$, we recall that the $\alpha$-Toeplitz matrices of
dimension $n\times n$ are defined as
\begin{eqnarray}\label{tnalpha}
  T_{n,\alpha}=[a_{r-\alpha c}]_{r,c=0}^{n-1},
\end{eqnarray}
where the quantities $r-\alpha s$ are not reduced modulus $n$. In
analogy with the case of $\alpha=1$, the elements $a_{j}$ are the
Fourier coefficients of some function $f$ in $L^{1}(Q)$, with
$Q=(-\pi,\pi)$, i.e., $a_j=\tilde f_j$ as in (\ref{defcoeff}) with
$d=1$.
%\begin{eqnarray}\label{bj}
%  b_{j}=\frac{1}{2\pi}\int_{-\pi}^{\pi}f(t)e^{-ijt}dt.
%\end{eqnarray}
If we denote by $T_{n}$ the classical Toeplitz matrix generated by
the function $f\in L^{1}(Q)$, $T_{n}=[a_{r-c}]_{r,c=0}^{n-1}$,
$a_{j}=\tilde f_j$ defined as in $(\ref{defcoeff})$, and by
$T_{n,\alpha}$ the $\alpha$-Toeplitz matrix generated by the same
function, one verifies immediately for $n$ and $\alpha$ generic that
\begin{eqnarray}\label{T}
  T_{n,\alpha}=\left[\widehat{T}_{n,\alpha}|\mathcal{T}_{n,\alpha}\right]=
  \left[T_{n}\widehat{Z}_{n,\alpha}|\mathcal{T}_{n,\alpha}\right],
\end{eqnarray}
where $\widehat{T}_{n,\alpha}\in\mathbb{C}^{n\times \mu_{\alpha}}$,
$\mu_{\alpha}=\left\lceil \frac{n}{\alpha}\right\rceil$, is the matrix $T_{n,\alpha}$
defined in $(\ref{tnalpha})$ by considering only the $\mu_{\alpha}$ first columns,
$\mathcal{T}_{n,\alpha}\in\mathbb{C}^{n\times (n-\mu_{\alpha})}$ is the matrix $T_{n,\alpha}$
defined in $(\ref{tnalpha})$ by considering only the $n-\mu_{\alpha}$ last columns,
and $\widehat{Z}_{n,\alpha}$ is the matrix defined in $(\ref{i})$ by considering
only the $\mu_{\alpha}$ first columns.
\begin{proof} (of relation $(\ref{T}).$) For $r=0,1,\ldots,n-1$ and $s=0,1,\ldots,\mu_{\alpha}-1,$ one has
\begin{eqnarray*}
    (\widehat{T}_{n,\alpha})_{r,s}&=&(T_{n})_{r,\alpha s},\\
    (\widehat{Z}_{n,\alpha})_{r,s}&=&\delta_{r-\alpha s},
\end{eqnarray*}
and
\begin{eqnarray*}
  (T_{n}\widehat{Z}_{n,\alpha})_{r,s} &=& \overset{n-1}
  {\underset{l=0}\sum}(T_{n})_{r,l}(\widehat{Z}_{n,\alpha})_{l,s} \\
    &=& \overset{n-1}{\underset{l=0}\sum}\delta_{l-\alpha s}(T_{n})_{r,l}  \\
    &{\underset{\rm (a)}=}& (T_{n})_{r,\alpha s} \\
    &=& (\widehat{T}_{n,\alpha})_{r,s},
\end{eqnarray*}
where (a) follows because there exists a unique
$l\in\{0,1,\ldots,n-1\}$ such that $l-\alpha s\equiv0\textrm{ (mod
$n$)}$, that is, $l\equiv \alpha s\textrm{ (mod $n$)}$, and, since
$0\leq\alpha s\leq n-1$, we obtain $l=\alpha s$.
\end{proof}

If we take the matrix $\widehat{T}_{n,\alpha}$ of size
$n\times(\mu_{\alpha}+1)$, then relation $(\ref{T})$ is no longer
true. In reality, looking at the $(\mu_{\alpha}+1)$-th column of the
$\alpha $-Toeplitz we observe Fourier coefficients with indices
which are not present (less or equal to $-n$) in the Toeplitz matrix
$T_{n}$. More precisely,
\begin{eqnarray*}
  (T_{n,\alpha})_{0,\mu_{\alpha}}=a_{0-\alpha\mu_{\alpha}}=
  a_{-\alpha\mu_{\alpha}},\qquad\textrm{and $\,-\alpha\mu_{\alpha}\leq -n$.}
\end{eqnarray*}

It follows that $\mu_{\alpha}$ is the maximum number of columns for
which relation $(\ref{T})$ is true.

\subsection{Some preparatory results}

We begin with some preliminary notations and definitions.

\begin{definition} \label{appr:seq}
Suppose a sequence of matrices $\{A_n\}_{n}$ of size $d_n$ is given.
We say that $\{\{B_{n,m}\}_{n}:m\ge 0\}$, $B_{n,m}$ of size $d_n$,
$m\in\mathbb{N}$, is an approximating class of sequences $(a.c.s.)$
for $\{A_n\}_{n}$ if, for all sufficiently large $m\in\mathbb{N}$,
the following splitting holds:
\begin{equation}\label{spli1}
A_n=B_{n,m}+R_{n,m}+N_{n,m}\quad \mbox{for all}\ n> n_m,
\end{equation}
with
\begin{equation}
\label{propri1} {\rm Rank}(R_{n,m}) \leq  d_n\, c(m), \quad \norma{
N_{n,m}} \leq  \omega(m),
\end{equation}
where $\|\cdot\|$ is the spectral norm (largest singular value),
$n_m$, $ c(m)$ and $\omega(m)$ depend only on $m$ and, moreover,
\begin{equation}
\label{propri1bis}
\lim_{m\to\infty} \omega(m)=0, \ \ \lim_{m\to\infty} c(m)=0.
\end{equation}
\end{definition}

%\begin{definition} \label{distribution}
%Suppose a sequence of matrices $\{A_n\}_n$ of size $d_n$ is given
%together with a measurable complex-valued function $\theta$. We say that $\{A_n\}_n$ is distributed in the singular value sense as the function $\theta$, and we write $\{A_n\}_n\sim_\sigma \theta$, if for
%every $F\in \Test(\mathbb{R}^{+}_{0})$ (continuous with bounded support)
%\[
%\lim_{n\to\infty} \frac 1 {d_n}\sum_{i=1}^{d_n}
%F\Bigl(\sigma_i\bigl(A_n\bigr)\Bigr)= \frac 1 {m(K)}\int_K
%F\bigl(|\theta(t)|\bigr)\, dt,
%\]
%with $\left\{\sigma_i\bigl(A_n\bigr)\right\}$ denoting the set of
%the singular values of $A_n$ (counted with their multiplicities) and
%with $K\subset {\mathbb C}^d$, $d\geq 1$,  being the definition set
%of $\theta$, with positive and finite Lebesgue measure $m(K)$.
%\end{definition}

\begin{proposition}{\rm \cite{algebra}}
\label{lem} Let $\{d_n\}_n$ be an increasing sequence of natural numbers.
Suppose a sequence formed by matrices $\{A_n\}_n$ of size $d_n$ is
given such that $\{\{B_{n,m}\}_n:\ m\ge 0\}$, $m\in\hat {\mathbb{N}}
\subset \mathbb{N}$, $\# \hat {\mathbb{N}}=\infty$, is an $a.c.s.$
for $\{A_n\}_n$ in the sense of Definition $\ref{appr:seq}$. Suppose
that $\{B_{n,m}\}_n\sim_\sigma (\theta_m, G)$ and that $\theta_m$
converges in measure to the measurable function $\theta$ over $G$.
Then necessarily
\begin{equation}
\label{l7} \{A_n\}_n\sim_\sigma (\theta,G),
\end{equation}
(see Definition $\ref{def-distribution}$).
\end{proposition}

\begin{proposition}\label{fgbis}{\rm \cite{algebra,ser-glt}}
If $\{A_{n}\}_{n}$ and $\{B_{n}\}_{n}$ are two sequences of matrices
of strictly increasing dimension, such that
$\{A_{n}\}_{n}\sim_{\sigma}(\theta,G)$ and
$\{B_{n}\}_{n}\sim_{\sigma} (0,G)$, then
\begin{eqnarray*}
  \{A_{n}+B_{n}\}_{n}\sim_{\sigma} (\theta,G).
\end{eqnarray*}
\end{proposition}

%\begin{theorem}
%Let $f$ be an $L^{1}$ function, that is a complex valued integrable function, over $Q=[-\pi,\pi)$, if $\{T_{n}(f)\}_{n}$ is the sequence of Toeplitz matrices generated by $f$, then it holds
%\begin{eqnarray*}
%  \{T_{n}(f)\}_{n}\sim_{\sigma}(f,Q).
%\end{eqnarray*}
%\end{theorem}

\begin{proposition}\label{fg}{\rm \cite{algebra}}
Let $f,\,g\in L^{1}(Q^d)$, $Q=(-\pi,\pi)$, and let $\{T_{n}(f)\}_{n}$ and
$\{T_{n}(g)\}_{n}$ be the two sequences of Toeplitz matrices
generated by $f$ and $g$, respectively.  The following
distribution result is true
\begin{eqnarray*}
  \{T_{n}(f)T_{n}(g)\}_{n}\sim_{\sigma}(fg,Q^d).
\end{eqnarray*}
\end{proposition}

\begin{lemma}\label{terzo} Let $f$ be a measurable complex-valued function on a set $K$,
and consider the measurable function $\sqrt{|f|}:K\rightarrow\mathbb{R}^{+}$.
Let $\{A_{n,m}\}$, with $A_{n,m}\in\mathbb{C}^{d_{n}\times d'_{n}}$, $d'_{n}\leq d_{n}$,
be a sequence of matrices of strictly increasing dimension: $d'_{n}<d'_{n+1}$ and $d_{n}\leq d_{n+1}$.
If the sequence of matrices $\{A_{n,m}^{*}A_{n,m}\}$, with
 $A_{n,m}^{*}A_{n,m}\in\mathbb{C}^{d'_{n}\times d'_{n}}$ and $d'_{n}<d'_{n+1}$, is distributed
 in the singular value sense as the function $f$ over a proper set $G\subset K$ in the sense of Definition $\ref{def-distribution}$,
 then the sequence $\{A_{n,m}\}$ is distributed in the singular value sense as the function $\sqrt{|f|}$
 over the same $G$.
\end{lemma}

\begin{proof} From the singular value decomposition ($SVD$), we can write $A_{n,m}$ as
\begin{eqnarray*}
  A_{n,m}=U\Sigma V^{*}=U\left[\begin{array}{cccc}
  \sigma_{1} & & &\\
  & \sigma_{2} & & \\
  & & \ddots & \\
  & & & \sigma_{d'_{n}}\\
  \hline
  & 0 & &
  \end{array}\right]V^{*},
\end{eqnarray*}
with $U$ and $V$ unitary matrices $U\in\mathbb{C}^{d_{n}\times d_{n}}$,
$V\in\mathbb{C}^{d'_{n}\times d'_{n}}$ and $\Sigma\in\mathbb{R}^{d_{n}\times d'_{n}}$, $\sigma_{j}\geq 0$;
by multiplying $A_{n,m}^{*}A_{n,m}$ we obtain:
\begin{eqnarray}\label{asa}
  \notag A_{n,m}^{*}A_{n,m}=V\Sigma^{T} U^{*}U\Sigma V^{*}&=&V\Sigma^{T}\Sigma V^{*}=V\Sigma^{(2)} V^{*}\\
  &=&V\left[\begin{array}{cccc}
  \sigma_{1}^{2} & & &\\
  & \sigma_{2}^{2} & & \\
  & & \ddots & \\
  & & & \sigma_{d'_{n}}^{2}
  \end{array}\right]V^{*},
\end{eqnarray}
with $V$ unitary matrix $V\in\mathbb{C}^{d'_{n}\times d'_{n}}$ and
$\Sigma^{(2)}\in\mathbb{R}^{d'_{n}\times d'_{n}}$, $\sigma_{j}^{2}\geq
0$; we observe that $(\ref{asa})$ is an $SVD$ for
$A_{n,m}^{*}A_{n,m}$, that is, the singular values
$\sigma_{j}(A_{n,m}^{*}A_{n,m})$ of $A_{n,m}^{*}A_{n,m}$ are the
square of singular values $\sigma_{j}(A_{n,m})$ of $A_{n,m}$. Since
$\{A_{n,m}^{*}A_{n,m}\}\sim_{\sigma} (f,G)$, by definition it hold
that for every $F\in \Test(\mathbb{R}_{0}^{+})$
\begin{eqnarray}\label{disasa}
\notag\lim_{n\to\infty} \frac{1}{d'_{n}}\sum_{i=1}^{d'_{n}}
F\left(\sigma_{i}(A_{n,m}^{*}A_{n,m})\right)&=&
\frac{1}{\mu(G)}\int_{G}
F\left(|f(t)|\right)\, dt\\
&=& \frac{1}{\mu(G)}\int_{G}H\left(\sqrt{|f(t)|}\right)\, dt,
\end{eqnarray}
where $H$ is such that $F=H\circ\sqrt{\cdot}$; but, owing to
$\sigma_{j}(A_{n,m})=\sqrt{\sigma_{j}(A_{n,m}^{*}A_{n,m})}$ we
obtain
\begin{eqnarray}\label{disasa1}
\notag\lim_{n\to\infty} \frac{1}{d'_{n}}\sum_{i=1}^{d'_{n}}F\left(\sigma_{i}(A_{n,m}^{*}A_{n,m})\right)&=&
\lim_{n\to\infty} \frac{1}{d'_{n}}\sum_{i=1}^{d'_{n}}F\left(\sigma_{i}^{2}(A_{n,m})\right)\\
&=&\lim_{n\to\infty} \frac{1}{d'_{n}}\sum_{i=1}^{d'_{n}}
H\left(\sigma_{i}(A_{n,m})\right).
\end{eqnarray}

From $(\ref{disasa})$ and $(\ref{disasa1})$ we obtain
\begin{eqnarray}
\lim_{n\to\infty} \frac{1}{d'_{n}}\sum_{i=1}^{d'_{n}}
H\left(\sigma_{i}(A_{n,m})\right)=\frac{1}{\mu(G)}\int_{G}
H\left(\sqrt{|f(t)|}\right)\, dt,
\end{eqnarray}
for every $H\in \Test(\mathbb{R}_{0}^{+})$, so
$\{A_{n,m}\}\sim_{\sigma}(\sqrt{|f(t)|},G)$.
\end{proof}

\begin{lemma}\label{primo} Let $\{A_{n}\}_{n}$ and $\{Q_{n}\}_{n}$ be two sequences of matrices
of strictly increasing dimension
($A_{n},Q_{n}\in\mathbb{C}^{d_{n}\times d_{n}}$, $d_{n}<d_{n+1}$),
where $Q_{n}$ are all unitary matrices ($Q_{n}Q_{n}^{*}=I$). If
$\{A_{n}\}_{n}\sim_{\sigma}(0, G)$ then
$\{A_{n}Q_{n}\}_{n}\sim_{\sigma} (0,G)$ and
$\{Q_{n}A_{n}\}_{n}\sim_{\sigma} (0,G)$.
\end{lemma}

\begin{proof} Putting $B_{n}=A_{n}Q_{n}$,  assuming that
\begin{eqnarray*}
  A_{n}=U_{n}\Sigma_{n}V_{n},
\end{eqnarray*}
is an $SVD$ for $A_{n}$, and taking into account that the product of
two unitary matrices is still a unitary matrix, we deduce that the
writing
\begin{eqnarray*}
  B_{n}=A_{n}Q_{n}=U_{n}\Sigma_{n}V_{n}Q_{n}=U_{n}\Sigma_{n}\widehat{V}_{n},
\end{eqnarray*}
is an $SVD$ for $B_{n}$. The latter implies that $A_{n}$ and $B_{n}$
have exactly the same singular values, so that the two sequences
$\{A_{n}\}_{n}$ and $\{B_{n}\}_{n}$ are distributed in the same way.
\end{proof}

\begin{lemma}\label{secondo} Let $\{A_{n}\}_{n}$ and $\{Q_{n}\}_{n}$ be two sequences of matrices of strictly
increasing dimension  ($A_{n},Q_{n}\in\mathbb{C}^{d_{n}\times
d_{n}}$, $d_{n}<d_{n+1}$). If $\{A_{n}\}_{n}\sim_{\sigma} (0,G)$ and
$\|Q_{n}\| \le M$ for some nonnegative constant $M$ independent of
$n$, then $\{A_{n}Q_{n}\}_{n}\sim_{\sigma} (0,G)$ and
$\{Q_{n}A_{n}\}_{n}\sim_{\sigma} (0,G)$.
\end{lemma}

\begin{proof} Since $\{A_{n}\}_{n}\sim_{\sigma} (0,G)$, then $\{0_{n}\}_{n}$ (sequence of zero matrices)
is an $a.c.s.$ for $\{A_{n}\}_{n}$; this means (by Definition $(\ref{appr:seq})$)
that we can write, for every $m$ sufficiently large, $m\in\mathbb{N}$
\begin{eqnarray}\label{A}
  A_{n}=0_{n}+R_{n,m}+N_{n,m},\qquad \forall n>n_{m},
\end{eqnarray}
with
\begin{eqnarray*}
  {\rm Rank}(R_{n,m})\leq d_{n}c(m),\qquad\|N_{n,m}\|\leq\omega(m),
\end{eqnarray*}
where $n_{m}\geq 0$, $c(m)$ and $\omega(m)$ depend only on $m$ and, moreover
\begin{eqnarray*}
  \lim_{m\rightarrow\infty}c(m)=0,\qquad\lim_{m\rightarrow\infty}\omega(m)=0.
\end{eqnarray*}

Now consider the matrix $A_{n}Q_{n}$; from $(\ref{A})$ we obtain
\begin{eqnarray*}
  A_{n}Q_{n}=0_{n}+R_{n,m}Q_{n}+N_{n,m}Q_{n},\qquad \forall n>n_{m},
\end{eqnarray*}
with
\begin{eqnarray*}
  &&{\rm Rank}(R_{n,m}Q_{n})\leq\min\{{\rm Rank}(R_{n,m}),{\rm Rank}(Q_{n})\}\leq
  {\rm Rank}(R_{n,m})\leq d_{n}c(m),\\
  &&\|N_{n,m}Q_{n}\|\leq\|N_{n,m}\|\|Q_{n}\|\leq M\omega(m),
\end{eqnarray*}
where
\begin{eqnarray*}
  \lim_{m\rightarrow\infty}c(m)=0,\qquad\lim_{m\rightarrow\infty}M\omega(m)=0,
\end{eqnarray*}
then $\{0_{n}\}_{n}$ is an $a.c.s.$ for the sequence
$\{A_{n}Q_{n}\}_{n}$ and, by Proposition $\ref{lem}$,
$\{A_{n}Q_{n}\}_{n}\sim_{\sigma} (0,G)$.
\end{proof}

\subsection{Singular value distribution for the $\alpha$-Toeplitz sequences}

As stated in formula $(\ref{T})$, the matrix $T_{n,\alpha}$ can be
written as
\begin{eqnarray}\label{dist}
  \notag T_{n,\alpha}&=&\left[T_{n}\widehat{Z}_{n,\alpha}|\mathcal{T}_{n,\alpha}\right]\\
  &=&\left[\begin{array}{c|c}
  T_{n}\widehat{Z}_{n,\alpha} & 0
  \end{array}\right]+\left[\begin{array}{c|c}
  0 & \mathcal{T}_{n,\alpha}
  \end{array}\right].
\end{eqnarray}

To find the distribution in the singular value sense of the sequence
$\{T_{n,\alpha}\}_{n}$, the idea is to study separately the
distribution of the two sequences
$\{[T_{n}\widehat{Z}_{n,\alpha}|0]\}_{n}$ and
$\{[0|\mathcal{T}_{n,\alpha}]\}_{n}$, to prove
$\{[0|\mathcal{T}_{n,\alpha}]\}_{n}\sim (0,G)$, and then apply
Proposition $\ref{fgbis}$.

\subsubsection{Singular value distribution for the sequence $\{[T_{n}\widehat{Z}_{n,\alpha}|0]\}_{n}$}

Since $T_{n}\widehat{Z}_{n,\alpha}\in\mathbb{C}^{n\times
\mu_{\alpha}}$ and
$[T_{n}\widehat{Z}_{n,\alpha}|0]\in\mathbb{C}^{n\times n}$, the
matrix $[T_{n}\widehat{Z}_{n,\alpha}|0]$ has $n-\mu_{\alpha}$
singular values equal to zero and the remaining $\mu_{\alpha}$ equal
to those of $T_{n}\widehat{Z}_{n,\alpha}$; to study the distribution
in the singular value sense of this sequence of non-square matrices,
we use Lemma $\ref{terzo}$: consider the $\alpha$-Toeplitz matrix
``truncated''
$\widehat{T}_{n,\alpha}=T_{n}(f)\widehat{Z}_{n,\alpha}$, where the
elements of the Toeplitz matrix $T_{n}(f)=[a_{r-c}]_{r,c=0}^{n-1}$
are the Fourier coefficients of a function $f$ in $L^{1}(Q)$, $Q=(-\pi,\pi)$, then
we have
\begin{eqnarray}\label{ff}
  \notag\widehat{T}_{n,\alpha}^{*}\widehat{T}_{n,\alpha}&=&(T_{n}(f)\widehat{Z}_{n,\alpha})^{*}T_{n}(f)\widehat{Z}_{n,\alpha}=
  \widehat{Z}_{n,\alpha}^{*}T_{n}(f)^{*}T_{n}(f)\widehat{Z}_{n,\alpha}\\
  &=&\widehat{Z}_{n,\alpha}^{*}T_{n}(\overline{f})T_{n}(f)\widehat{Z}_{n,\alpha}.
\end{eqnarray}

We provide in detail the analysis in the case where  $f\in
L^{2}(Q)$. The general setting in which $f\in L^{1}(Q)$ can be
obtained by approximation and density arguments as done in
\cite{algebra}. From Proposition $\ref{fg}$ if $f\in L^{2}(Q)\subset
L^{1}(Q)$ (that is $|f|^{2}\in L^{1}(Q)$), then
$\{T_{n}(\overline{f})T_{n}(f)\}_{n}\sim_{\sigma}(|f|^{2},Q)$.
Consequently, for every $m$ sufficiently large, $m\in\mathbb{N}$,
the use of Proposition $\ref{lem}$ implies
\begin{eqnarray*}
  T_{n}(\overline{f})T_{n}(f)=T_{n}(|f|^{2})+R_{n,m}+N_{n,m},\qquad \forall n>n_{m},
\end{eqnarray*}
with
\begin{eqnarray*}
  {\rm Rank}(R_{n,m})\leq nc(m),\qquad\|N_{n,m}\|\leq\omega(m),
\end{eqnarray*}
where $n_{m}\geq 0$, $c(m)$ and $\omega(m)$ depend only on $m$ and, moreover
\begin{eqnarray*}
  \lim_{m\rightarrow\infty}c(m)=0,\qquad\lim_{m\rightarrow\infty}\omega(m)=0.
\end{eqnarray*}

Therefore $(\ref{ff})$ becomes
\begin{eqnarray}\label{tst}
  \notag\widehat{T}_{n,\alpha}^{*}\widehat{T}_{n,\alpha}
  &=&\widehat{Z}_{n,\alpha}^{*}(T_{n}(|f|^{2})+R_{n,m}+N_{n,m})\widehat{Z}_{n,\alpha}\\
  \notag&=&\widehat{Z}_{n,\alpha}^{*}T_{n}(|f|^{2})\widehat{Z}_{n,\alpha}+
  \widehat{Z}_{n,\alpha}^{*}R_{n,m}\widehat{Z}_{n,\alpha}+\widehat{Z}_{n,\alpha}^{*}N_{n,m}\widehat{Z}_{n,\alpha}\\
  &=&\widehat{Z}_{n,\alpha}^{*}T_{n}(|f|^{2})\widehat{Z}_{n,\alpha}+\widehat{R}_{n,m,\alpha}+\widehat{N}_{n,m,\alpha},
\end{eqnarray}
with
\begin{eqnarray}\label{1bis}
  &&{\rm Rank}(\widehat{R}_{n,m,\alpha})\leq\min\{{\rm Rank}(\breve{Z}_{n,\alpha}),{\rm Rank}(R_{n,m})\}\leq
  {\rm Rank}(R_{n,m})\leq nc(m),\\
  \label{2bis}&&\|\widehat{N}_{n,m,\alpha}\|\leq2\|\breve{Z}_{n,\alpha}\|\|N_{n,m}\|\leq2\omega(m),
\end{eqnarray}
and
\begin{eqnarray*}
  \lim_{m\rightarrow\infty}c(m)=0,\qquad\lim_{m\rightarrow\infty}2\omega(m)=0,
\end{eqnarray*}
where in $(\ref{1bis})$ and $(\ref{2bis})$,
$\breve{Z}_{n,\alpha}=[\widehat{Z}_{n,\alpha}|0]\in\mathbb{C}^{n\times
n}$. In other words $\breve{Z}_{n,\alpha}$ is the matrix
$\widehat{Z}_{n,\alpha}$ supplemented by an appropriate number of
zero columns in order to make it square. Furthermore, it is worth
noticing that
$\|\widehat{Z}_{n,\alpha}\|=\|\widehat{Z}_{n,\alpha}^*\|=1$, because
$\widehat{Z}_{n,\alpha}$ is a submatrix of the identity: we have
used the latter relations in $(\ref{2bis})$.

Now, consider the matrix
$\widehat{Z}_{n,\alpha}^{*}T_{n}(|f|^{2})\widehat{Z}_{n,\alpha}
\in\mathbb{C}^{\mu_{\alpha}\times\mu_{\alpha}}$, with
$\mu_{\alpha}=\left\lceil \frac{n}{\alpha}\right\rceil$, $f\in
L^{2}(Q)\subset L^{1}(Q)$ (so $|f|^{2}\in L^{1}(Q)$). From
$(\ref{T})$, setting
$T_{n}=T_{n}(|f|^{2})=[\tilde{a}_{r-c}]_{r,c=0}^{n-1}$, with
$\tilde{a}_{j}$ being the Fourier coefficients of $|f|^{2}$, and
setting $T_{n,\alpha}$ the $\alpha$-Toeplitz generated by the same
function $|f|^{2}$, it is immediate to observe
\begin{eqnarray}\label{brc}
  T_{n}\widehat{Z}_{n,\alpha}=\widehat{T}_{n,\alpha}\in
  \mathbb{C}^{n\times \mu_{\alpha}},\qquad \textrm{with}
  \quad (\widehat{T}_{n,\alpha})_{r,c}=\tilde{a}_{r-\alpha c},
\end{eqnarray}
for $r=0,\ldots,n-1$ and $c=0,\ldots,\mu_{\alpha}-1$. If we compute
$\widehat{Z}_{n,\alpha}^{*}\widehat{T}_{n,\alpha}\in\mathbb{C}^{\mu_{\alpha}\times\mu_{\alpha}}$,
where $Z_{n,\alpha}^{*}=[\delta_{c-\alpha r}]_{r,c=0}^{n-1}$ ($\delta_{k}$ defined as in $(\ref{i})$)
and $\widehat{Z}_{n,\alpha}^{*}\in\mathbb{C}^{\mu_{\alpha}\times n}$ is the submatrix of $Z_{n,\alpha}^{*}$ obtained by considering only the $\mu_{\alpha}$ first rows, for $r,c=0,\ldots,\mu_{\alpha}-1$, we obtain
\begin{eqnarray*}
  (\widehat{Z}_{n,\alpha}^{*}T_{n}(|f|^{2})\widehat{Z}_{n,\alpha})_{r,c}&=&
  (\widehat{Z}_{n,\alpha}^{*}\widehat{T}_{n,\alpha})_{r,c}\\
  &=&\sum_{\ell=0}^{n-1}(\widehat{Z}_{n,\alpha}^{*})_{r,\ell}(\widehat{T}_{n,\alpha})_{\ell,c}\\
  &{\underset{\rm (a)}=}&(\widehat{T}_{n,\alpha})_{\alpha r,c}\\
  &{\underset{\rm from\,(\ref{brc})}=}& \widehat{a}_{\alpha r-\alpha c},
\end{eqnarray*}
where (a) follows from the existence of a unique
$\ell\in\{0,1,\ldots,n-1\}$ such that $\ell-\alpha
r\equiv0\textrm{ (mod $n$)}$, that is, $\ell\equiv \alpha
r\textrm{ (mod $n$)}$, and, since $0\leq\alpha r\leq n-1$, we find
$\ell=\alpha r$.

Therefore
\begin{eqnarray*}
  \widehat{Z}_{n,\alpha}^{*}T_{n}(|f|^{2})\widehat{Z}_{n,\alpha}&=&
  [\tilde{a}_{\alpha r-\alpha c}]_{r,c=0}^{\mu_{\alpha}-1}\\
  &=&T_{\mu_{\alpha}}(\widehat{|f|^{(2)}}),
\end{eqnarray*}
where $\widehat{|f|^{(2)}}\in L^{1}(Q)$ is given by
\begin{eqnarray}\label{f2t}
  \widehat{|f|^{(2)}}(x)&=&\frac{1}{\alpha}\sum_{j=0}^{\alpha-1}|f|^{2}\left(\frac{x+2\pi j}{\alpha}\right),\\
  \label{f2tbis}|f|^{2}(x)&=&\sum_{k=-\infty}^{+\infty}\tilde{a}_{k}e^{ikx}.
\end{eqnarray}

\begin{proof} (of relation $(\ref{f2t}).$) We denote by $\textrm{a}_{j}$ the Fourier coefficients of
$\widehat{|f|^{(2)}}$. We want to show that for
$r,c=0,\ldots,\mu_{\alpha}-1$, $\textrm{a}_{r-c}=\tilde{a}_{\alpha
r-\alpha c}$, where $\tilde{a}_{k}$ are the Fourier coefficients of
$|f|^{2}$. From $(\ref{defcoeff})$, $(\ref{f2t})$ and
$(\ref{f2tbis})$, we have
\begin{eqnarray*}
  \textrm{a}_{r-c}&=&\frac{1}{2\pi}\int_{-\pi}^{\pi}\frac{1}{\alpha}\sum_{j=0}^{\alpha-1}
  \sum_{k=-\infty}^{+\infty}\tilde{a}_{k}e^{ik\left(\frac{x+2\pi j}{\alpha}\right)}e^{-i(r-c)x}dx\\
  &=&\frac{1}{2\pi\alpha}\int_{-\pi}^{\pi}\sum_{k=-\infty}^{+\infty}\tilde{a}_{k}\left(\sum_{j=0}^{\alpha-1}e^{\frac{i2\pi kj}{\alpha}}\right)
  e^{\frac{ikx}{\alpha}}e^{-i(r-c)x}dx.
\end{eqnarray*}
Some remarks are in order:
\begin{itemize}
  \item[-] if $k$ is a multiple of $\alpha$, $k=\alpha t$ for some value of $t$, then we have that $\overset{\alpha-1}{\underset{j=0}\sum}e^{\frac{i2\pi kj}{\alpha}}=\overset{\alpha-1}{\underset{j=0}\sum}e^{\frac{i2\pi\alpha tj}{\alpha}}=
    \overset{\alpha-1}{\underset{j=0}\sum}e^{i2\pi tj}=\overset{\alpha-1}{\underset{j=0}\sum}1=\alpha$.
  \item[-] if $k$ is not a multiple of $\alpha$, then $e^{\frac{i2\pi k}{\alpha}}\neq 1$ and therefore $\overset{\alpha-1}{\underset{j=0}\sum}e^{\frac{i2\pi kj}{\alpha}}=\overset{\alpha-1}{\underset{j=0}\sum}\left(e^{\frac{i2\pi k}{\alpha}}\right)^{j}$ is a finite geometric series whose sum is given by
  \begin{eqnarray*}
    \sum_{j=0}^{\alpha-1}\left(e^{\frac{i2\pi k}{\alpha}}\right)^{j}=\frac{1-e^{\frac{i2\pi k\alpha}{\alpha}}}{1-e^{\frac{i2\pi k}{\alpha}}}=
    \frac{1-e^{i2\pi k}}{1-e^{\frac{i2\pi k}{\alpha}}}=\frac{1-1}{1-e^{\frac{i2\pi k}{\alpha}}}=0.
  \end{eqnarray*}
\end{itemize}

Finally, taking into account the latter statements and recalling
that $\frac{1}{2\pi}\int_{-\pi}^{\pi} e^{i\ell x} dx=
\left\{\begin{smallmatrix} 1 & \textrm{ if $\ell=0$}\\ 0
&\textrm{otherwise} \end{smallmatrix}\right.$, we find
\begin{eqnarray*}
  \textrm{a}_{r-c}&=&\frac{1}{2\pi\alpha}\int_{-\pi}^{\pi}\sum_{t=-\infty}^{+\infty}\tilde{a}_{\alpha t}
  \alpha e^{\frac{i\alpha t x}{\alpha}}e^{-i(r-c)x}dx\\
  &=&\sum_{t=-\infty}^{+\infty}\tilde{a}_{\alpha t}\frac{1}{2\pi}\int_{-\pi}^{\pi}e^{ix(t-(r-c))}dx\\
  &=&\tilde{a}_{\alpha(r-c)}.
\end{eqnarray*}
\end{proof}

In summary, from $(\ref{tst})$ we have
\begin{eqnarray*}
  \widehat{T}_{n,\alpha}^{*}\widehat{T}_{n,\alpha}=
  T_{\mu_{\alpha}}(\widehat{|f|^{(2)}})+\widehat{R}_{n,m,\alpha}+\widehat{N}_{n,m,\alpha},
\end{eqnarray*}
with
$\{T_{\mu_{\alpha}}(\widehat{|f|^{(2)}})\}_{n}\sim_{\sigma}(\widehat{|f|^{(2)}},Q)$.
We recall that, owing to $(\ref{f2t})$, the relation $|f|^{2}\in
L^{1}(Q)$ implies $\widehat{|f|^{(2)}}\in L^{1}(Q)$. Consequently
Proposition $\ref{lem}$ implies that
$\{\widehat{T}_{n,\alpha}^{*}\widehat{T}_{n,\alpha}\}_{n}\sim_{\sigma}(\widehat{|f|^{(2)}},Q)$.
Clearly $\widehat{|f|^{(2)}}\in L^{1}(Q)$ is equivalent to write
$\sqrt{\widehat{|f|^{(2)}}}\in L^{2}(Q)$: therefore, from Lemma
$\ref{terzo}$, we infer
$\{\widehat{T}_{n,\alpha}\}_{n}\sim_{\sigma}(\sqrt{\widehat{|f|^{(2)}}},Q)$.

Now, as mentioned at the beginning of this section, by Definition
$\ref{def-distribution}$, we have
\begin{eqnarray*}
\lim_{n\rightarrow\infty} \frac{1}{n}\sum_{j=1}^{n}F\left(\sigma_j([\widehat{T}_{n,\alpha}|0])\right)&=&
\lim_{n\rightarrow\infty}\frac{1}{n}\sum_{j=1}^{\mu_{\alpha}}F\left(\sigma_j([\widehat{T}_{n,\alpha}|0])\right)+
\lim_{n\rightarrow\infty}\frac{1}{n}\sum_{j=\mu_{\alpha}+1}^{n}F(0)\\
&=&\lim_{n\rightarrow\infty}\frac{\mu_{\alpha}}{n}\sum_{j=1}^{\mu_{\alpha}}\frac{F\left(\sigma_j([\widehat{T}_{n,\alpha}|0])\right)}{\mu_{\alpha}}+
\lim_{n\rightarrow\infty}\frac{n-\mu_{\alpha}}{n}F(0)\\
&=&\frac{1}{\alpha}\frac{1}{2\pi}\int_{-\pi}^{\pi}F\left(\sqrt{\widehat{|f|^{(2)}}(x)}\right)\,dx+
\left(1-\frac{1}{\alpha}\right)F(0),
\end{eqnarray*}
which results to be equivalent to the following distribution formula
\begin{eqnarray}\label{dist1}
  \{[T_{n}\widehat{Z}_{n,\alpha}|0]\}_{n}\sim_{\sigma}(\theta, Q\times[0,1]),
\end{eqnarray}
where
\begin{eqnarray}\label{teta}
  \theta(x,t)=\left\{\begin{array}{cl}
  \sqrt{\widehat{|f|^{(2)}}(x)} & \textrm{for $t\in\left[0,\frac{1}{\alpha}\right]$,}\\
  0 & \textrm{for $t\in\left(\frac{1}{\alpha},1\right]$}.
  \end{array}\right.
\end{eqnarray}

\subsubsection{Singular value distribution for the sequence $\{[0|\mathcal{T}_{n,\alpha}]\}_{n}$}

In perfect analogy with the case of the matrix
$[T_{n}\widehat{Z}_{n,\alpha}|0]$, we can observe that
$\mathcal{T}_{n,\alpha}\in\mathbb{C}^{n\times (n-\mu_{\alpha})}$ and
$[0|\mathcal{T}_{n,\alpha}]\in\mathbb{C}^{n\times n}$. Therefore the
matrix $[0|\mathcal{T}_{n,\alpha}]$ has $\mu_{\alpha}$ singular
values equal to zero and the remaining $n-\mu_{\alpha}$ equal to
those of $\mathcal{T}_{n,\alpha}$. However, in this case we have
additional difficulties with respect to the matrix
$\widehat{T}_{n,\alpha}=T_{n}\widehat{Z}_{n,\alpha}$, because it is
not always true that $\mathcal{T}_{n,\alpha}$ can be written as
$T_{n}\mathcal{Z}_{n,\alpha}$, where $\mathcal{Z}_{n,\alpha}$ is the
matrix obtained by considering the $n-\mu_{\alpha}$ last columns of
$Z_{n,\alpha}$. Indeed, in $\mathcal{T}_{n,\alpha}$ there are
Fourier coefficients with index, in modulus, greater than $n$: the
Toeplitz matrix $T_{n}=[a_{r-c}]_{r,c=0}^{n-1}$ has coefficients
$a_j$ with $j$ ranging from $1-n$ to $n-1$, while the
$\alpha$-Toeplitz matrix $T_{n,\alpha}=[a_{r-\alpha
c}]_{r,c=0}^{n-1}$ has $a_{n-1}$ as coefficient of maximum index and
$a_{-\alpha(n-1)}$ as coefficient of minimum index, and, if
$\alpha\geq 2$, we have $-\alpha(n-1)<-(n-1)$.

Even if we take the Toeplitz matrix $T_{n}$, which has as its first
column the first column of $\mathcal{T}_{n,\alpha}$ and the other
generated according to the rule $(T_{n})_{j,k}=a_{j-k}$, it is not
always true that we can write $\mathcal{T}_{n,\alpha}=T_{n}P$ for a
suitable submatrix $P$ of a permutation matrix, indeed, if the matrix
$T_{n}=[\beta_{r-c}]_{r,c=0}^{n-1}$ has as first column the first
column of $\mathcal{T}_{n,\alpha}$, we find that
$\beta_{0}=(\mathcal{T}_{n,\alpha})_{0,0}=(T_{n,\alpha})_{0,\mu_{\alpha}}=a_{-\alpha\mu_{\alpha}}$.
As a consequence, $T_{n}$ has
$\beta_{-(n-1)}=a_{-(n-1)-\alpha\mu_{\alpha}}$ as coefficient of
minimum index, while $\mathcal{T}_{n,\alpha}$ has $a_{-\alpha(n-1)}$
as coefficient of minimum index. Therefore
\begin{eqnarray*}
  -(n-1)\alpha-(-(n-1)-\alpha\mu_{\alpha})&=&(1-\alpha)(n-1)+\alpha\mu_{\alpha}\qquad\;\;\; n\leq\alpha\mu_{\alpha}=\alpha\left\lceil \frac{n}{\alpha}\right\rceil\leq (n+\alpha-1)\\
  &\leq &(1-\alpha)(n-1)+(n+\alpha-1)\\
  &=&(1-\alpha)(n-1)+(n-1)+\alpha\\
  &=&(n-1)(1-\alpha+1)+\alpha\\
  &=&(2-\alpha)(n-1)+\alpha< 0\qquad\; \textrm{for $\alpha>2$ and
  $n>4$}.
\end{eqnarray*}

Thus, if $\alpha>2$ and $n>4$ we have
$-(n-1)\alpha<-(n-1)-\alpha\mu_{\alpha}$ and the coefficient of
minimum index $a_{-\alpha(n-1)}$ of $\mathcal{T}_{n,\alpha}$ is not
contained in the matrix $T_{n}$ that has
$a_{-(n-1)-\alpha\mu_{\alpha}}$ as coefficient of minimum index.

Then we proceed in another way: in the first column of $\mathcal{T}_{n,\alpha}\in\mathbb{C}^{n\times (n-\mu_{\alpha})}$
(and consequently throughout the matrix) there are only coefficients with index $< 0$,
indeed coefficient with the largest index of $\mathcal{T}_{n,\alpha}$ is
$(\mathcal{T}_{n,\alpha})_{n-1,0}=(T_{n,\alpha})_{n-1,\mu_{\alpha}}=a_{n-1-\alpha\mu_{\alpha}}$ and
$n-1-\alpha\mu_{\alpha}\leq n-1-n<0$ and the coefficient with smallest index is
$(\mathcal{T}_{n,\alpha})_{0,n-\mu_{\alpha}-1}=(T_{n,\alpha})_{0,n-\mu_{\alpha}-1+\mu_{\alpha}}=
(T_{n,\alpha})_{0,n-1}=a_{-\alpha(n-1)}$. Consider therefore a Toeplitz matrix $T_{d_{n,\alpha}}$
of dimension $d_{n,\alpha}$ with $d_{n,\alpha}>\frac{\alpha(n-1)}{2}+1$, defined in this way:
\begin{eqnarray}\label{tdn}
   T_{d_{n,\alpha}}=\left[\begin{array}{ccccc}
     a_{-d_{n,\alpha}+1} & a_{-d_{n,\alpha}} & a_{-d_{n,\alpha}-1} & \cdots & a_{-2d_{n,\alpha}+2}   \\
     a_{-d_{n,\alpha}+2} & a_{-d_{n,\alpha}+1} & \ddots & \ddots & a_{-2d_{n,\alpha}+3} \\
     \vdots & \ddots & \ddots & \ddots & \vdots\\
     a_{-1} & a_{-2} & \ddots & \ddots & a_{-d_{n,\alpha}}\\
     a_{0} & a_{-1} & a_{-2} & \cdots &  a_{-d_{n,\alpha}+1}
   \end{array}\right]=\left[a_{r-c-d_{n,\alpha}+1}\right]_{r,c=0}^{d_{n,\alpha}-1}.
\end{eqnarray}

Since the coefficient with smallest index is $a_{-2d_{n,\alpha}+2}$,
we find
\begin{eqnarray*}
  -2d_{n,\alpha}+2<-2\left(\frac{\alpha(n-1)}{2}+1\right)+2=-\alpha(n-1)-2+2=-\alpha(n-1).
\end{eqnarray*}

As a consequence, we obtain that all the coefficients of
$\mathcal{T}_{n,\alpha}$ are ``contained'' in the matrix
$T_{d_{n,\alpha}}$. In particular, if
\begin{eqnarray*}
d_{n,\alpha}>(\alpha-1)(n-1)+2,
\end{eqnarray*}
(this condition ensures $d_{n,\alpha}>\frac{\alpha(n-1)}{2}+1$, that
all the subsequent inequalities are correct, and that the size of
all the matrices involved are non-negative), then it can be shown
that
\begin{eqnarray}\label{hat}
  \mathcal{T}_{n,\alpha}=\left[0_{1}|I_{n}|0_{2}\right]T_{d_{n,\alpha}}\mathcal{Z}_{d_{n,\alpha},\alpha},
\end{eqnarray}
where
$\mathcal{Z}_{d_{n,\alpha},\alpha}\in\mathbb{C}^{d_{n,\alpha}\times
(n-\mu_{\alpha})}$ is the matrix defined in (\ref{i}), of dimension
$d_{n,\alpha}\times d_{n,\alpha}$, by considering only the
$n-\mu_{\alpha}$ first columns and
$\left[0_{1}|I_{n}|0_{2}\right]\in\mathbb{C}^{n\times d_{n,\alpha}}$
is a block matrix with $0_{1}\in\mathbb{C}^{n\times
(d_{n,\alpha}-\alpha\mu_{\alpha}-1)}$ and
$0_{2}\in\mathbb{C}^{n\times (\alpha\mu_{\alpha}-n+1)}$.
%$I_{n}\in\mathbb{C}^{n\times n}$ is the identity matrix.

\begin{proof} (of relation $(\ref{hat}).$) First we observe that:
\begin{itemize}
  \item[] for $r=0,1,\ldots,n-1$ and $s=0,1,\ldots,n-\mu_{\alpha}-1$ we have
  \begin{eqnarray}\label{eq1}
    (\mathcal{T}_{n,\alpha})_{r,s}=(T_{n,\alpha})_{r,s+\mu_{\alpha}}=a_{r-\alpha
    s-\alpha\mu_{\alpha}};
  \end{eqnarray}
  \item[] for $r=0,1,\ldots,n-1$ and $s=0,1,\ldots,d_{n,\alpha}-1$ we have
  \begin{eqnarray}\label{eq2}
    (\left[0_{1}|I_{n}|0_{2}\right])_{r,s}=\left\{\begin{array}{cl} 1 &
    \textrm{if $s=r+d_{n,\alpha}-\alpha\mu_{\alpha}-1$},\\
    0 & \textrm{otherwise};\end{array}\right.
  \end{eqnarray}
  \item[] for $r,s=0,1,\ldots,d_{n,\alpha}-1$ we have
  \begin{eqnarray*}
    (T_{d_{n,\alpha}})_{r,s}=a_{r-s-d_{n,\alpha}+1};
  \end{eqnarray*}
  \item[] for $r=0,1,\ldots,d_{n,\alpha}-1$ and $s=0,1,\ldots,n-\mu_{\alpha}-1,$ we have
  \begin{eqnarray*}
    (\mathcal{Z}_{d_{n,\alpha},\alpha})_{r,s}=\delta_{r-\alpha s}.
  \end{eqnarray*}
\end{itemize}

Since
$T_{d_{n,\alpha}}\mathcal{Z}_{d_{n,\alpha},\alpha}\in\mathbb{C}^{d_{n,\alpha}\times
(n-\mu_{\alpha})}$, for $r=0,1,\ldots,d_{n,\alpha}-1$ and
$s=0,1,\ldots,n-\mu_{\alpha}-1,$ it holds
\begin{eqnarray}\label{eq3}
  \notag(T_{d_{n,\alpha}}\mathcal{Z}_{d_{n,\alpha},\alpha})_{r,s} &=& \overset{d_{n,\alpha}-1}
  {\underset{l=0}\sum}(T_{d_{n,\alpha}})_{r,l}(\mathcal{Z}_{d_{n,\alpha},\alpha})_{l,s} \\
  \notag&=& \overset{d_{n,\alpha}-1}{\underset{l=0}\sum}\delta_{l-\alpha s}a_{r-l-d_{n,\alpha}+1}\\
  & {\underset{\rm (a)}=} & a_{r-\alpha s-d_{n,\alpha}+1},
\end{eqnarray}
where (a) follows from the existence of a unique
$l\in\{0,1,\ldots,d_{n,\alpha}-1\}$ such that $l-\alpha
s\equiv 0\textrm{ (mod $d_{n,\alpha}$)}$, that is, $l\equiv \alpha
s\textrm{ (mod $d_{n,\alpha}$)}$, and, since $0\leq\alpha s\leq d_{n,\alpha}-1$, we have $l=\alpha s$. Since
$\left[0_{1}|I_{n}|0_{2}\right]T_{d_{n,\alpha}}\mathcal{Z}_{d_{n,\alpha},\alpha}\in\mathbb{C}^{n\times
(n-\mu_{\alpha})}$, for $r=0,1,\ldots,n-1$ and
$s=0,1,\ldots,n-\mu_{\alpha}-1,$ we find
\begin{eqnarray*}
  (\left[0_{1}|I_{n}|0_{2}\right]T_{d_{n,\alpha}}\mathcal{Z}_{d_{n,\alpha},\alpha})_{r,s} &=&
  \overset{d_{n,\alpha}-1}{\underset{l=0}\sum}(\left[0_{1}|I_{n}|0_{2}\right])_{r,l}
  (T_{d_{n,\alpha}}\mathcal{Z}_{d_{n,\alpha},\alpha})_{l,s} \\
  & {\underset{\rm (d)}=}&a_{r+d_{n,\alpha}-\alpha\mu_{\alpha}-1-\alpha s-d_{n,\alpha}+1}\\
  &=&a_{r-\alpha\mu_{\alpha}-\alpha s}\\
  &{\underset{\rm from\,(\ref{eq1})}=}&(\mathcal{T}_{n,\alpha})_{r,s},
\end{eqnarray*}
where (d) follows from $(\ref{eq3})$,
$(T_{d_{n,\alpha}}\mathcal{Z}_{d_{n,\alpha},\alpha})_{l,s}=a_{l-\alpha
s-d_{n,\alpha}+1}$, and from the following fact: using
$(\ref{eq2})$, we find $(\left[0_{1}|I_{n}|0_{2}\right])_{r,l}=1$ if
and only if $l=r+d_{n,\alpha}-\alpha\mu_{\alpha}-1$.
\end{proof}

We can now observe immediately that the matrix $T_{d_{n,\alpha}}$ defined in $(\ref{tdn})$ can be written as
\begin{eqnarray}\label{flip}
  T_{d_{n,\alpha}}=JH_{d_{n,\alpha}},
\end{eqnarray}
where $J$ is the ``flip'' matrix of dimension $d_{n,\alpha}\times d_{n,\alpha}$:
\begin{eqnarray*}
  J=\left[\begin{array}{cccc}
   &  &  & 1\\
   & & 1 & \\
   & \adots & & \\
   1 &  & &
  \end{array}\right],
\end{eqnarray*}
and $H_{d_{n,\alpha}}$ is the Hankel matrix of dimension $d_{n,\alpha}\times d_{n,\alpha}$:
\begin{eqnarray*}
   H_{d_{n,\alpha}}=\left[\begin{array}{ccccc}
     a_{0} & a_{-1} & a_{-2} & \cdots &  a_{-d_{n,\alpha}+1}  \\
     a_{-1} & a_{-2} & \adots & \adots & a_{-d_{n,\alpha}}\\
     \vdots & \adots & \adots & \adots & \vdots\\
     a_{-d_{n,\alpha}+2} & a_{-d_{n,\alpha}+1} & \adots & \adots & a_{-2d_{n,\alpha}+3} \\
     a_{-d_{n,\alpha}+1} & a_{-d_{n,\alpha}} & a_{-d_{n,\alpha}-1} & \cdots & a_{-2d_{n,\alpha}+2}
   \end{array}\right].
\end{eqnarray*}

If $f(x)\in L^{1}(Q)$, $Q=(-\pi,\pi)$, is the generating function of the Toeplitz matrix $T_{n}=T_{n}(f)=[a_{r-c}]_{r,c=0}^{n-1}$ in $(\ref{T})$, where the $k$-th Fourier coefficient of $f$ is $a_k$,
%
%\begin{eqnarray*}
%  f(x)=\sum_{k=-\infty}^{+\infty}a_{k}e^{ikx},
%\end{eqnarray*}
%
then $f(-x)\in L^{1}(Q)$ is the generating function of the Hankel
matrix $H_{d_{n,\alpha}}=[a_{-r-c}]_{r,c=0}^{d_{n,\alpha}-1}$; by
invoking Theorem 6, page 161 of \cite{FasTi}, the sequence of
matrices $\{H_{d_{n,\alpha}}\}$ is distributed in the singular value
sense as the zero function: $\{H_{d_{n,\alpha}}\}\sim_{\sigma}
(0,Q)$. From Lemma $\ref{primo}$, by $(\ref{flip})$, since $J$ is a
unitary matrix, we have $\{T_{d_{n,\alpha}}\}\sim_{\sigma} (0,Q)$ as
well.

Consider the decomposition in $(\ref{hat})$:
\begin{eqnarray*}
  \mathcal{T}_{n,\alpha}=\left[0_{1}|I_{n}|0_{2}\right]T_{d_{n,\alpha}}\mathcal{Z}_{d_{n,\alpha},\alpha}=
  Q_{d_{n,\alpha}}T_{d_{n,\alpha}}\mathcal{Z}_{d_{n,\alpha},\alpha}.
\end{eqnarray*}

If we complete the matrices $Q_{d_{n,\alpha}}\in\mathbb{C}^{n\times
d_{n,\alpha}}$ and
$\mathcal{Z}_{d_{n,\alpha},\alpha}\in\mathbb{C}^{d_{n,\alpha}\times
(n-\mu_{\alpha})}$ by adding an appropriate number of zero rows and
columns, respectively, in order to make it square
\begin{eqnarray*}
  \mathbf{Q}_{d_{n,\alpha}}&=&\left[\begin{array}{c}
  Q_{d_{n,\alpha}}\\
  \hline
  0
  \end{array}\right]\in\mathbb{C}^{d_{n,\alpha}\times d_{n,\alpha}},\\
  \mathbf{Z}_{d_{n,\alpha},\alpha}&=&\left[\begin{array}{c|c}
  \mathcal{Z}_{d_{n,\alpha},\alpha} & 0
  \end{array}\right]\in\mathbb{C}^{d_{n,\alpha}\times d_{n,\alpha}},
\end{eqnarray*}
then it is immediate to note that
\begin{eqnarray*}
  \mathbf{Q}_{d_{n,\alpha}}T_{d_{n,\alpha}}\mathbf{Z}_{d_{n,\alpha},\alpha}=\left[\begin{array}{c|c}
  \mathcal{T}_{n,\alpha} & 0 \\
  \hline
  0 & 0
  \end{array}\right]=\mathbf{T}_{n,\alpha}\in\mathbb{C}^{d_{n,\alpha}\times d_{n,\alpha}}.
\end{eqnarray*}

From Lemma $\ref{secondo}$, since
$\|\mathbf{Q}_{d_{n,\alpha}}\|=\|\mathbf{Z}_{d_{n,\alpha},\alpha}\|=1$
(indeed they are both ``incomplete'' permutation matrices), and
since $\{T_{d_{n,\alpha}}\}\sim_{\sigma} (0,Q)$,  we infer that
$\{\mathbf{T}_{n,\alpha}\}\sim_{\sigma} (0,Q)$.

Recall that $\mathbf{T}_{n,\alpha}\in\mathbb{C}^{d_{n,\alpha}\times
d_{n,\alpha}}$ with $d_{n,\alpha}>(\alpha-1)(n-1)+2$; then we
can always choose $d_{n,\alpha}$ such that $\alpha
n=d_{n,\alpha}>(\alpha-1)(n-1)+2$ (if $n,\alpha\geq2$). Now, since
$\{\mathbf{T}_{n,\alpha}\}\sim_{\sigma} (0,Q)$, it holds that the
sequence $\{\mathbf{T}_{n,\alpha}\}$ is weakly clustered at zero in
the singular value sense, i.e., $\forall\epsilon>0$,
\begin{eqnarray}\label{sigmabft}
  \sharp\{j:\sigma_{j}(\mathbf{T}_{n,\alpha})>\epsilon\}=o(d_{n,\alpha})=o(\alpha n)=o(n).
\end{eqnarray}

The matrix $\mathbf{T}_{n,\alpha}$ is a block matrix that can be
written as
\begin{eqnarray*}
  \mathbf{T}_{n,\alpha}=\left[\begin{array}{c|c}
  \mathcal{T}_{n,\alpha} & 0 \\
  \hline
  0 & 0
  \end{array}\right]=\left[\begin{array}{c|c}
  [\mathcal{T}_{n,\alpha}|0] & 0 \\
  \hline
  0 & 0
  \end{array}\right],
\end{eqnarray*}
where $\mathcal{T}_{n,\alpha}\in\mathbb{C}^{n\times (n-\mu_{\alpha})}$
and $[\mathcal{T}_{n,\alpha}|0]\in\mathbb{C}^{n\times n}$. By the
singular value decomposition we obtain
\begin{eqnarray*}
  \mathbf{T}_{n,\alpha}=\left[\begin{array}{c|c}
  [\mathcal{T}_{n,\alpha}|0] & 0 \\
  \hline
  0 & 0
  \end{array}\right]=\left[\begin{array}{c|c}
  U_{1}\Sigma_{1} V_{1}^{*} & 0 \\
  \hline
  0 & U_{2}0V_{2}^{*}
  \end{array}\right]=\left[\begin{array}{c|c}
  U_{1} & 0 \\
  \hline
  0 & U_{2}
  \end{array}\right]\left[\begin{array}{c|c}
  \Sigma_{1} & 0 \\
  \hline
  0 & 0
  \end{array}\right]\left[\begin{array}{c|c}
  V_{1} & 0 \\
  \hline
  0 & V_{2}
  \end{array}\right]^{*},
\end{eqnarray*}
that is, the singular values of $\mathbf{T}_{n,\alpha}$ that are
different from zero are the singular values of
$[\mathcal{T}_{n,\alpha}|0]\in\mathbb{C}^{n\times n}$. Thus
$(\ref{sigmabft})$ can be written as follows: $\forall\epsilon>0$,
\begin{eqnarray*}
  \sharp\{j:\sigma_{j}([\mathcal{T}_{n,\alpha}|0])>\epsilon\}=o(d_{n,\alpha})=o(\alpha
  n)=o(n).
\end{eqnarray*}

The latter relation means that the sequence
$\{[\mathcal{T}_{n,\alpha}|0]\}_{n}$ is weakly clustered at zero in
the singular value sense, and hence
$\{[\mathcal{T}_{n,\alpha}|0]\}_{n}\sim_{\sigma} (0,Q)$. If we now
consider the matrix
\begin{eqnarray*}
  \hat G=\left[\begin{array}{c|c}
  0 & I_{n-\mu_{\alpha}}\\
  \hline
  0 & 0
  \end{array}\right]\in\mathbb{C}^{n\times n},
\end{eqnarray*}
where $I_{n-\mu_{\alpha}}$ is the identity matrix of dimension
$(n-\mu_{\alpha})\times(n-\mu_{\alpha})$, then
$[\mathcal{T}_{n,\alpha}|0]\hat G=[0|\mathcal{T}_{n,\alpha}]$, and
since $\|\hat G\|=1$ and
$\{[\mathcal{T}_{n,\alpha}|0]\}_{n}\sim_{\sigma} (0,Q)$, from Lemma
$\ref{secondo}$ we find
\begin{eqnarray}\label{dist2}
\{[0|\mathcal{T}_{n,\alpha}]\}_{n}\sim_{\sigma} (0,Q).
\end{eqnarray}

In conclusion: from the relations $(\ref{dist})$, $(\ref{dist1})$
and $(\ref{dist2})$, using Proposition $\ref{fgbis}$ with $G=Q\times
[0,1]$, we obtain that
\begin{eqnarray*}
  \{T_{n,\alpha}\}_{n}\sim_{\sigma}(\theta, Q\times [0,1]),
\end{eqnarray*}
where $\theta$ is defined in $(\ref{teta})$. Notice that for
$\alpha=1$ the symbol $\theta(x,t)$ coincides with $|f|(x)$ on the
extended domain $Q\times [0,1]$. Hence the
Szeg\"o-Tilli-Tyrtyshnikov-Zamarashkin result is found as a
particular case. Indeed $\theta(x,t)=|f|(x)$ does not depend on $t$
and therefore this additional variable can be suppressed i.e.
$\{T_{n,\alpha}\}_{n}\sim_{\sigma}(f,Q)$ with $T_{n,\alpha}=T_n(f)$.
The fact that the distribution formula is not unique should not
surprise since this phenomenon is inherent to the measure theory
because any measure-preserving exchange function is a distribution
function if one representative of the class is.

\section{Some remarks on multigrid methods}\label{sec:multigrid}

In the design of multigrid methods for large positive definite
linear systems one of the key points is to maintain the structure
(if any) of the original matrix in the lower levels. This means that
at every recursion level the new projected linear system should
retain the main properties of the original matrix (e.g. bandedness,
the same level of conditioning, the same algebra/Toeplitz/graph
structure etc.). Here for the sake of simplicity the example that
has to be considered is the one-level circulant case. Following
\cite{ADS,mcirco}, if $A_n=C_n$ is a positive circulant matrix of
size $n$ with $n$ power of $2$, then the projected matrix $A_k$ with
$k=n/2$ is defined as
\begin{equation}\label{proj}
A_k=\widetilde{Z}_{n,2}^T P_n^* A_n P_n \widetilde{Z}_{n,2},
\end{equation}
where $P_n$ is an additional circulant matrix. It is worth noticing that
the structure is kept since for every circulant $P_n$ the matrix
$A_k$ is a circulant matrix of size $k=n/2$. The features of the
specific $P_n$ have to be designed in such a way that the
convergence speed of the related multigrid is as high as possible
(see \cite{FS2,ADS} for a general strategy). We observe that the
eigenvalues of $A_k$ are given by
\begin{equation}\label{eig-proj}
\frac{1}{2}\overset{1}
  {\underset{l=0}\sum} g\left(\frac{x_j+2\pi l}{2} \right),
  \quad x_j=\frac{2\pi j}{k}, \ \ j=0,1,\ldots,k-1,\ k=n/2,
\end{equation}
where $g$ is the polynomial associated with the circulant matrix
$P_n^* A_n P_n$ in the sense of Subsection \ref{special-case}.
Therefore the singular values of $(P_n^* A_n
P_n)^{1/2}\widetilde{Z}_{n,2}$ are given by
\begin{equation}\label{sv-sqrt-proj}
\frac{1}{\sqrt{2}}\sqrt{\overset{1}
  {\underset{l=0}\sum} g\left(\frac{x_j+2\pi l}{2} \right)},
  \quad x_j=\frac{2\pi j}{k}, \ \ j=0,1,\ldots,k-1,\ k=n/2.
\end{equation}

Notice that the latter formula is a special instance of
(\ref{sv-nontrivial-symbol-p}) for $|p|^2=g$ ($g$ is necessarily
nonnegative since it can be written a $|q|^2 f$ where $q$ is the
polynomial associated with $P_n$ and $f$ the nonnegative polynomial
associated with $A_n$), for $\alpha=2$ and $n$ even number so that
$(n,2)=2$. Therefore, according to (\ref{sv-nontrivial-symbol-p}),
the numbers in (\ref{sv-sqrt-proj}) identify the nontrivial singular
values of the $2$-circulant matrix $(P_n^* A_n P_n)^{1/2}Z_{n,2}$ up to a scaling factor.
In other words $\alpha$-circulant matrices arise naturally in the
design of fast multigrid solvers for circulant linear systems and,
along the same lines, $\alpha$-Toeplitz matrices arise naturally in
the design of fast multigrid solvers for Toeplitz linear systems;
see \cite{FS2,ADS, Sun}.

Conversely, we now can see clearly that formula
(\ref{sv-nontrivial-symbol-p}) furnishes a wide generalization of
the spectral analysis of the projected matrices, by allowing a
higher degree of freedom: we can choose $n$ divisible by $\alpha$
with $\alpha\neq 2$, we can choose $n$ not divisible by $\alpha$.
Such a degree of freedom is not just academic, but could be
exploited for devising optimally convergent multigrid solvers also
in critical cases emphasized e.g. in \cite{ADS, Sun}. In particular,
if $x_0$ is an isolated zero of $f$ (the nonnegative polynomial
related to $A_n=C_n$) and also $\pi+x_0$ is a zero for the same
function, then due to special symmetries, the associated multigrid
(or even two-grid) method cannot be optimal. In other words, for
reaching a preassigned accuracy, we cannot expect a number of
iterations independent of the order $n$. However these pathological
symmetries are due to the choice of $\alpha=2$, so that a choice of a
projector as $P_n\widetilde{Z}_{n,\alpha}$ for a different
$\alpha\neq 2$ and a different $n$ could completely overcome the
latter drawback.

\section{Generalizations}\label{sec:gen}

First of all we observe that the
requirement that the symbol $f$ is square integrable can be removed.
In \cite{algebra} it is proven that the singular value distribution
of $\{T_{n}(f)T_n(g)\}_{n}$ is given by $h=fg$ with $f,g$ being just
Lebesgue integrable and with $h$ that is only measurable and
therefore may fail to be Lebesgue integrable. This fact is
sufficient for extending the proof of the relation
$\{T_{n,\alpha}\}_{n}\sim_{\sigma}(\theta,Q\times [0,1])$ to the
case where $\theta(x,t)$ is defined as in $(\ref{teta})$ with the original
symbol $f\in L^1$.

Now we consider the general multilevel case. When $\alpha$ is a positive vector, we have
\begin{eqnarray}\label{dist1-d}
  \{T_{n,\alpha}\}_{n}\sim_{\sigma}(\theta,Q^d\times[0,1]^d),
\end{eqnarray}
where
\begin{eqnarray}\label{teta-d}
  \theta(x,t)=\left\{\begin{array}{cl}
  \sqrt{\widehat{|f|^{(2)}}(x)} & \textrm{for $t\in\left[\underline{0},\frac{1}{\alpha}\right]$,}\\
  0 & \textrm{for $t\in\left(\frac{1}{\alpha},e\right]$},
  \end{array}\right.
\end{eqnarray}
with
\begin{eqnarray}\label{f2t-d}
  \widehat{|f|^{(2)}}(x)&=&\frac{1}{\hat{\alpha}}\sum_{j=\underline{0}}^{\alpha-e}|f|^{2}\left(\frac{x+2\pi j}{\alpha}\right),
%  \notag|f|^{2}(x)&=&\sum_{k=-\infty}^{+\infty}\widehat{b}_{k}e^{ikx}.
\end{eqnarray}
and where all the arguments are modulus $2\pi$ and all the
operations are intended componentwise that is
$t\in\left[\underline{0},\frac{1}{\alpha}\right]$ means that $t_k\in
[0,1/\alpha_k]$, $k=1,\ldots,d$,
$t\in\left(\frac{1}{\alpha},e\right]$ means that $t_k\in
(1/\alpha_k,1]$, $k=1,\ldots,d$, the writing $\frac{x+2\pi
j}{\alpha}$ defines the $d$-dimensional vector whose $k$-th
component is $(x_j+2\pi j_k)/\alpha_k$, $k=1,\ldots,d$, and $\hat{\alpha}=\alpha_{1}\alpha_{2}\cdots\alpha_{d}$.

\subsubsection{Examples of $\alpha$-circulant and $\alpha$-Toeplitz matrices when some of the entries of $\alpha$ vanish}

We start this subsection with a brief digression on multilevel matrices.
A $d$-level matrix $A$ of dimension $\hat{n}\times\hat{n}$ with $n=(n_{1},n_{2},\ldots,n_{d})$ and $\hat{n}=n_{1}n_{2}\cdots n_{d}$ can be viewed as a matrix of dimension $n_{1}\times n_{1}$ in which each element is a block of dimension $n_{2}n_{3}\cdots n_{d}\times n_{2}n_{3}\cdots n_{d}$; in turn, each block of dimension $n_{2}n_{3}\cdots n_{d}\times n_{2}n_{3}\cdots n_{d}$ can be viewed as a matrix of dimension $n_{2}\times n_{2}$ in which each element is a block of dimension $n_{3}n_{4}\cdots n_{d}\times n_{3}n_{4}\cdots n_{d}$, and so on. So we can say that $n_{1}$ is the most ``outer'' dimension of the matrix $A$ and $n_{d}$ is the most ``inner'' dimension. If we multiply by an appropriate permutation matrix $P$ the $d$-level matrix $A$, we can exchange the ``order of dimensions'' of $A$,
namely $P^{T}AP$ becomes a matrix again of dimension
$\hat{n}\times\hat{n}$ but with $n=(n_{p(1)},n_{p(2)},\ldots,n_{p(d)})$ and $\hat{n}=n_{p(1)}n_{p(2)}\cdots n_{p(d)}=n_{1}n_{2}\cdots n_{d}$ (where $p$ is a permutation of $d$ elements) and $n_{p(1)}$ is the most ``outer'' dimension of the matrix $A$ and $n_{p(d)}$ is the most ``inner'' dimension.

This trick helps us to understand what happens to the singular values of $\alpha$-circulant and $\alpha$-Toeplitz $d$-level matrices, especially when some of the entries of the vector $\alpha$ are zero; indeed, as we observed in Subsection $\ref{alphazero}$, if $\alpha=\underline{0}$, the $d$-level $\alpha$-circulant (or $\alpha$-Toeplitz) matrix $A$ is a block matrix with constant blocks on each row, so if we order the vector $\alpha$ (which has some components equal to zero) so that the components equal to zero are in the top positions, $\alpha=(0,\ldots,0,\alpha_{k},\ldots,\alpha_{d})$, the matrix $P^{T}AP$ (where $P$ is the permutation matrix associated with $p$) becomes a block matrix with constant blocks on each row and with blocks of dimension $n_{k}\cdots n_{d}\times n_{k}\cdots n_{d}$; with this ``new'' structure, formulas $(\ref{eq-2-1-3})$ and $(\ref{eq-2-1-3-bis})$ are even more intuitively understandable, as we shall see later in the examples.
\begin{lemma}\label{scambio} Let $A$ be a 2-level Toeplitz matrix of dimension $\hat{n}\times\hat{n}$ with $n=(n_{1},n_{2})$ and
$\hat{n}=n_{1}n_{2}$,
\begin{eqnarray*}
  A=\left[\left[a_{(j_{1}-k_{1},j_{2}-k_{2})}\right]_{j_{2},k_{2}=0}^{n_{2}-1}\right]_{j_{1},k_{1}=0}^{n_{1}-1}.
\end{eqnarray*}
There exists a permutation matrix $P$ such that
\begin{eqnarray*}
  P^{T}AP=\left[\left[a_{(j_{1}-k_{1},j_{2}-k_{2})}\right]_{j_{1},k_{1}=0}^{n_{1}-1}\right]_{j_{2},k_{2}=0}^{n_{2}-1}.
\end{eqnarray*}
\end{lemma}

\begin{example} Let $n=(n_{1},n_{2})=(2,3)$ and consider the 2-level Toeplitz matrix $A$ of dimension $6\times6$
\begin{eqnarray*}
  A=\left[
  \begin{array}{ccc|ccc}
    a_{(0,0)} & a_{(0,-1)} & a_{(0,-2)} & a_{(-1,0)} & a_{(-1,-1)} & a_{(-1,-2)}\\
    a_{(0,1)} & a_{(0,0)} & a_{(0,-1)} & a_{(-1,1)} & a_{(-1,0)} & a_{(-1,-1)}\\
    a_{(0,2)} & a_{(0,1)} & a_{(0,0)} & a_{(-1,2)} & a_{(-1,1)} & a_{(-1,0)}\\
    \hline
    a_{(1,0)} & a_{(1,-1)} & a_{(1,-2)} & a_{(0,0)} & a_{(0,-1)} & a_{(0,-2)}\\
    a_{(1,1)} & a_{(1,0)} & a_{(1,-1)} & a_{(0,1)} & a_{(0,0)} & a_{(0,-1)}\\
    a_{(1,2)} & a_{(1,1)} & a_{(1,0)} & a_{(0,2)} & a_{(0,1)} & a_{(0,0)}
  \end{array}\right].
\end{eqnarray*}

This matrix can be viewed as a matrix of dimension $2\times 2$ in which each element is a block of dimension $3\times 3$. If we take the permutation matrix
\begin{eqnarray*}
  P=\left[\begin{array}{cccccc}
  1 & 0 & 0 & 0 & 0 & 0 \\
  0 & 0 & 1 & 0 & 0 & 0 \\
  0 & 0 & 0 & 0 & 1 & 0 \\
  0 & 1 & 0 & 0 & 0 & 0 \\
  0 & 0 & 0 & 1 & 0 & 0 \\
  0 & 0 & 0 & 0 & 0 & 1
  \end{array}\right],
\end{eqnarray*}
then it is plain to see that
\begin{eqnarray*}
  P^{T}AP=\left[
  \begin{array}{cc|cc|cc}
    a_{(0,0)} & a_{(-1,0)} & a_{(0,-1)} & a_{(-1,-1)} & a_{(0,-2)} & a_{(-1,-2)}\\
    a_{(1,0)} & a_{(0,0)}  & a_{(1,-1)} & a_{(0,-1)}  & a_{(1,-2)} & a_{(0,-2)}\\
    \hline
    a_{(0,1)} & a_{(-1,1)} & a_{(0,0)} & a_{(-1,0)} & a_{(0,-1)} & a_{(-1,-1)}\\
    a_{(1,1)} & a_{(0,1)}  & a_{(1,0)} & a_{(0,0)}  & a_{(1,-1)} & a_{(0,-1)}\\
    \hline
    a_{(0,2)} & a_{(-1,2)} & a_{(0,1)} & a_{(-1,1)} & a_{(0,0)} & a_{(-1,0)}\\
    a_{(1,2)} & a_{(0,2)}  & a_{(1,1)} & a_{(0,1)}  & a_{(1,0)} & a_{(0,0)}
  \end{array}\right],
\end{eqnarray*}
and now $P^{T}AP$ can be naturally viewed as a matrix of dimension $3\times 3$ in which each element is a block of dimension $2\times 2 $.
\end{example}

\begin{corollary}\label{scambiod} Let $A$ be a $d$-level Toeplitz matrix of dimension $\hat{n}\times\hat{n}$ with $n=(n_{1},n_{2},\ldots,n_{d})$ and
$\hat{n}=n_{1}n_{2}\cdots n_{d}$,
\begin{eqnarray*}
A=\left[\left[\cdots\left[a_{(j_{1}-k_{1},j_{2}-k_{2},\ldots,j_{d}-k_{d})}\right]_{j_{d},k_{d}=0}^{n_{d}-1}\cdots\right]_{j_{2},k_{2}=0}^{n_{2}-1}
\right]_{j_{1},k_{1}=0}^{n_{1}-1}.
\end{eqnarray*}
For every permutation $p$ of $d$ elements, there exists a permutation matrix $P$ such that
\begin{eqnarray*}
  P^{T}AP=\left[\left[\cdots\left[a_{(j_{1}-k_{1},j_{2}-k_{2},\ldots,j_{d}-k_{d})}\right]_{j_{p(d)},k_{p(d)}=0}^{n_{p(d)}-1}\cdots\right]_{j_{p(2)},k_{p(2)}=0}^{n_{p(2)}-1}\right]_{j_{p(1)},k_{p(1)}=0}^{n_{p(1)}-1}.
\end{eqnarray*}
\end{corollary}

\begin{oss} Lemma $\ref{scambio}$ and Corollary $\ref{scambiod}$ also apply to $d$-level $\alpha$-circulant and $\alpha$-Toeplitz matrices.
\end{oss}

Now, let $\alpha=(\alpha_{1},\alpha_{2},\ldots,\alpha_{d})$ be a $d$-dimensional vector of nonnegative integers and $t=\sharp\{j:\alpha_{j}=0\}$ be the number of zero entries of $\alpha$. If we take a permutation $p$ of $d$ elements such that $\alpha_{p(1)}=\alpha_{p(2)}=\ldots=\alpha_{p(t)}=0$, (that is, $p$ is a permutation that moves all the zero components of the vector $\alpha$ in the top positions), then it is easy to prove that formulas $(\ref{eq-2-1-3})$ and $(\ref{eq-2-1-3-bis})$ remain the same for the matrix $P^{T}AP$ (where $P$ is the permutation matrix associated with $p$) but with $n[0]=(n_{p(1)},n_{p(2)},\ldots,n_{p(t)})$ and where $C_j$ and $T_j$ are a $d^+$-level $\alpha^+$-circulant and $\alpha^+$-Toeplitz matrix, respectively, with $\alpha^+=(\alpha_{p(t+1)},\alpha_{p(t+2)},\ldots,\alpha_{p(d)})$, of partial sizes $n[>0]=(n_{p(t+1)},n_{p(t+2)},\ldots,n_{p(d)})$, and whose expressions are
\begin{eqnarray*}
C_j&=&\left[\left[\cdots\left[a_{(r-\alpha \circ s)\ {\rm mod}\,
n}\right]_{r_{p(d)},s_{p(d)}=0}^{n_{p(d)}-1}\cdots\right]_{r_{p(t+2)},s_{p(t+2)}=0}^{n_{p(t+2)}-1}\right]_{r_{p(t+1)},s_{p(t+1)}=0}^{n_{p(t+1)}-1},\\
T_j&=&\left[\left[\cdots\left[a_{(r-\alpha \circ s)}\right]_{r_{p(d)},s_{p(d)}=0}^{n_{p(d)}-1}\cdots\right]_{r_{p(t+2)},s_{p(t+2)}=0}^{n_{p(t+2)}-1}\right]_{r_{p(t+1)},s_{p(t+1)}=0}^{n_{p(t+1)}-1},
\end{eqnarray*}
with $(r_{p(1)},r_{p(2)},\ldots,r_{p(t)})=j$. Obviously ${\rm
Sgval}(A)={\rm Sgval}(P^{T}AP)$.

We recall that if $B$ is a matrix of size $n \times n$ positive
semidefinite, that is $B^{*}=B$ and $x^{*}Bx\geq 0$ $\forall x\neq
0$, then ${\rm Eig}(B)={\rm Sgval}(B)$. Moreover, if $B=U\Sigma
U^{*}$ is a $SVD$ for $B$ (which coincides with the Schur
decomposition of $B$) with $\Sigma=\diag(\sigma_{j})$, then
\begin{eqnarray}\label{bmezzi}
  B^{1/2}=U\Sigma^{1/2}U^{*},
\end{eqnarray}
where $\Sigma^{1/2}=\diag(\sqrt{\sigma_{j}})$.

We proceed with two detailed examples: a 3-level $\alpha$-circulant matrix with $\alpha=(\alpha_{1},\alpha_{2},\alpha_{3})=(1,2,0)$, and a 3-level $\alpha$-Toeplitz with $\alpha=(\alpha_{1},\alpha_{2},\alpha_{3})=(0,1,2)$, which helps us to understand what happens if the vector $\alpha$ is not strictly positive. Finally we will propose the explicit calculation of the singular values of a $d$-level $\alpha$-circulant matrix in the particular case where the vector $\alpha$ has only one component different from zero.

\begin{example} Consider a 3-level $\alpha$-circulant matrix $A$ where $\alpha=(\alpha_{1},\alpha_{2},\alpha_{3})=(1,2,0)$
\begin{eqnarray*}
  A&=&\left[\left[\left[a_{((r_{1}-1\cdot s_{1})\textrm{ mod $n_{1}$},(r_{2}-2\cdot s_{2})\textrm{ mod $n_{2}$},(r_{3}-0\cdot s_{3})\textrm{ mod $n_{3}$})}
      \right]_{r_{3},s_{3}=0}^{n_{3}-1}\right]_{r_{2},s_{2}=0}^{n_{2}-1}
      \right]_{r_{1},s_{1}=0}^{n_{1}-1} \\
   &=&\left[\left[\left[a_{((r_{1}-s_{1})\textrm{ mod $n_{1}$},(r_{2}-2s_{2})\textrm{ mod $n_{2}$},r_{3})}
      \right]_{r_{3}=0}^{n_{3}-1}\right]_{r_{2},s_{2}=0}^{n_{2}-1}
      \right]_{r_{1},s_{1}=0}^{n_{1}-1}.
\end{eqnarray*}

If we choose a permutation $p$ of 3 elements such that
\begin{eqnarray*}
  &&(p(1),p(2),p(3))=(3,2,1),\\
  &&(\alpha_{p(1)},\alpha_{p(2)},\alpha_{p(3)})=(0,2,1), \\
  &&(n_{p(1)},n_{p(2)},n_{p(3)})=(n_{3},n_{2},n_{1}),
\end{eqnarray*}
and if we take the permutation matrix $P$ related to $p$, then
\begin{eqnarray*}
  P^{T}AP\equiv\hat{A}=\left[\left[\left[a_{((r_{1}-s_{1})\textrm{ mod $n_{1}$},(r_{2}-2s_{2})\textrm{ mod $n_{2}$},r_{3})}
      \right]_{r_{1},s_{1}=0}^{n_{1}-1}\right]_{r_{2},s_{2}=0}^{n_{2}-1}
      \right]_{r_{3}=0}^{n_{3}-1}.
\end{eqnarray*}

Now, for $r_{3}=0,1,...,n_{3}-1$, let us set
\begin{eqnarray*}
    C_{r_{3}}=\left[\left[a_{((r_{1}-s_{1})\textrm{ mod $n_{1}$},(r_{2}-2s_{2})\textrm{ mod $n_{2}$},r_{3})}
      \right]_{r_{1},s_{1}=0}^{n_{1}-1}\right]_{r_{2},s_{2}=0}^{n_{2}-1}.
\end{eqnarray*}
As a consequence, $C_{r_{3}}$ is a 2-level $\alpha^+$-circulant matrix with $\alpha^+=(2,1)$ and of partial
sizes $n[>0]=(n_{2},n_{1})$ and the matrix $\hat{A}$ can be rewritten as
\begin{eqnarray*}
  \hat{A}=\left[\begin{array}{cccc}
       C_{0} & C_{0} & \cdots & C_{0} \\
       C_{1} & C_{1} & \cdots & C_{1} \\
       \vdots & \vdots & \vdots & \vdots \\
       C_{n_{3}-1} & C_{n_{3}-1} & \cdots & C_{n_{3}-1}
    \end{array}\right],
\end{eqnarray*}
and this is a block matrix with constant blocks on each row.
From formula $(\ref{1})$, the singular values of $\hat{A}$ are the square
root of the eigenvalues of $\hat{A}^{*}\hat{A}$:
\begin{eqnarray*}
  \hat{A}^{*}\hat{A} &=&\left[\begin{array}{cccc}
                C_{0}^{\ast} & C_{1}^{\ast} & \cdots & C_{n_{3}-1}^{\ast} \\
                C_{0}^{\ast} & C_{1}^{\ast} & \cdots & C_{n_{3}-1}^{\ast} \\
                \vdots & \vdots & \vdots & \vdots \\
                C_{0}^{\ast} & C_{1}^{\ast} & \cdots & C_{n_{3}-1}^{\ast} \\
                \end{array}\right]\left[\begin{array}{cccc}
                C_{0} & C_{0} & \cdots & C_{0} \\
                C_{1} & C_{1} & \cdots & C_{1} \\
                \vdots & \vdots & \vdots & \vdots \\
                C_{n_{3}-1} & C_{n_{3}-1} & \cdots & C_{n_{3}-1}
                \end{array}\right] \\
    &=&\left[\begin{array}{cccc}
          \overset{n_{3}-1}{\underset{j=0}\sum}C_{j}^{\ast}C_{j} &
          \overset{n_{3}-1}{\underset{j=0}\sum}C_{j}^{\ast}C_{j} & \cdots &
          \overset{n_{3}-1}{\underset{j=0}\sum}C_{j}^{\ast}C_{j} \\
          \overset{n_{3}-1}{\underset{j=0}\sum}C_{j}^{\ast}C_{j} &
          \overset{n_{3}-1}{\underset{j=0}\sum}C_{j}^{\ast}C_{j} & \cdots &
          \overset{n_{3}-1}{\underset{j=0}\sum}C_{j}^{\ast}C_{j} \\
          \vdots & \vdots & \vdots & \vdots\\
          \overset{n_{3}-1}{\underset{j=0}\sum}C_{j}^{\ast}C_{j} &
          \overset{n_{3}-1}{\underset{j=0}\sum}C_{j}^{\ast}C_{j} & \cdots &
          \overset{n_{3}-1}{\underset{j=0}\sum}C_{j}^{\ast}C_{j} \\
        \end{array}\right] \\
   &=&  \underbrace{\left[\begin{array}{cccc}
         1 & 1 & \cdots& 1 \\
         1 & 1 & \cdots& 1 \\
         \vdots & \vdots & \vdots & \vdots \\
         1 & 1 &\cdots& 1 \\
       \end{array}\right]}_{\textrm{$n_{3}$ times}}\otimes\overset{n_{3}-1}{\underset{j=0}\sum}C_{j}^{\ast}C_{j}\\
   &=& J_{n_{3}}\otimes\overset{n_{3}-1}{\underset{j=0}\sum}C_{j}^{\ast}C_{j}.
\end{eqnarray*}

Therefore
\begin{equation}\label{AstarA}
    {\rm Eig}(\hat{A}^{*}\hat{A})={\rm Eig}\left(J_{n_{3}}\otimes\overset{n_{3}-1}{\underset{j=0}\sum}C_{j}^{\ast}C_{j}\right),
\end{equation}
where
\begin{eqnarray}\label{eigJ}
  {\rm Eig}(J_{n_{3}}) = \{0,n_{3}\},
\end{eqnarray}
because $J_{n_{3}}$ is a
matrix of rank 1, so it has all eigenvalues equal to zero except one
eigenvalue equal to ${\rm tr}(J_{n_{3}})=n_{3}$ (${\rm tr}$ is the trace of a matrix). If we put
\begin{eqnarray*}
  \lambda_{k}=\lambda_{k}\left(\overset{n_{3}-1}{\underset{j=0}\sum}C_{j}^{\ast}C_{j}\right),\qquad k=0,\ldots,n_{1}n_{2}-1,
\end{eqnarray*}
by exploiting basic properties of the tensor product and taking into
consideration $(\ref{AstarA})$ and $(\ref{eigJ})$ we find
\begin{eqnarray}\label{eigAA1}
  \lambda_{k}(\hat{A}^{*}\hat{A})&=&n_{3}\lambda_{k},\qquad k=0,\ldots,n_{1}n_{2}-1,\\
  \label{eigAA2}\lambda_{k}(\hat{A}^{*}\hat{A})&=&0,\qquad k=n_{1}n_{2},\ldots,n_{1}n_{2}n_{3}-1.
\end{eqnarray}

From $(\ref{eigAA1})$, $(\ref{eigAA2})$ and $(\ref{1}),$ and
recalling that ${\rm Sgval}(\hat{A})={\rm Sgval}(A)$, one obtains
that the singular values of $A$ are given by
\begin{eqnarray*}
  \sigma_{k}(A)&=&\sqrt{n_{3}\lambda_{k}},\qquad k=0,\ldots,n_{1}n_{2}-1,\\
  \sigma_{k}(A)&=&0,\qquad k=n_{1}n_{2},\ldots,n_{1}n_{2}n_{3}-1,
\end{eqnarray*}
and, since $\overset{n_{3}-1}{\underset{j=0}\sum}C_{j}^{\ast}C_{j}$ is a positive semidefinite matrix, from $(\ref{bmezzi})$ we can write
\begin{eqnarray*}
  \sigma_{k}(A)&=&\sqrt{n_{3}}\widetilde{\sigma}_{k},\qquad k=0,\ldots,n_{1}n_{2}-1,\\
  \sigma_{k}(A)&=&0,\qquad k=n_{1}n_{2},\ldots,n_{1}n_{2}n_{3}-1,
\end{eqnarray*}
where $\widetilde{\sigma}_{k}$ are the singular values of $\left(\overset{n_{3}-1}{\underset{j=0}\sum}C_{j}^{\ast}C_{j}\right)^{1/2}$.

Regarding the distribution in the sense of singular values, let $F\in C_{0}(\mathbb{R}_{0}^{+})$, continuous function over
$\mathbb{R}_{0}^{+}$ with bounded support, then there exists $a\in\mathbb{R}^{+}$ such that
\begin{equation}\label{222}
    \left|F(x)\right|\leq a \text{\,\,\,}\forall x\in\mathbb{R}_{0}^{+}.
\end{equation}

From formula $(\ref{sigmaFA})$ we have
\begin{eqnarray*}
  \Sigma_{\sigma}(F,A_{n}) &=& \frac{1}{n_{1}n_{2}n_{3}}\overset{n_{1}n_{2}n_{3}-1}{\underset{k=0}\sum}
  F(\sqrt{n_{3}}\widetilde{\sigma}_{k}) \\
   &=& \frac{n_{1}n_{2}(n_{3}-1)F(0)}{n_{1}n_{2}n_{3}}+ \frac{1}{n_{1}n_{2}n_{3}}
   \overset{n_{1}n_{2}-1}{\underset{k=0}\sum}F(\sqrt{n_{3}}\widetilde{\sigma}_{k}) \\
   &=& \left(1-\frac{1}{n_{3}}\right)F(0)+ \frac{1}{n_{1}n_{2}n_{3}}
   \overset{n_{1}n_{2}-1}{\underset{k=0}\sum}F(\sqrt{n_{3}}\widetilde{\sigma}_{k}).
\end{eqnarray*}

According to $(\ref{222}),$ we find
\begin{equation*}
    -an_{1}n_{2}\leq\overset{n_{1}n_{2}-1}{\underset{k=0}\sum}F(\sqrt{n_{3}}\widetilde{\sigma}_{k})\leq
    an_{1}n_{2}.
\end{equation*}
Therefore
\begin{equation*}
    -\frac{a}{n_{3}}\leq\frac{1}{n_{1}n_{2}n_{3}}
   \overset{n_{1}n_{2}-1}{\underset{k=0}\sum}F(\sqrt{n_{3}}\widetilde{\sigma}_{k})\leq\frac{a}{n_{3}},
\end{equation*}
so that
\begin{equation*}
    \left(1-\frac{1}{n_{3}}\right)F(0)-\frac{a}{n_{3}}\leq\Sigma_{\sigma}(F,A_{n})\leq\left(1-\frac{1}{n_{3}}\right)F(0)+
    \frac{a}{n_{3}}.
\end{equation*}

Now, recalling that the writing $n\rightarrow\infty$ means $\min_{1\leq j\leq 3}n_{j}\rightarrow\infty$, we obtain
\begin{equation*}
    F(0)\leq\underset{n\rightarrow\infty}{\lim}\Sigma_{\sigma}(F,A_{n})\leq F(0),
\end{equation*}
which implies
\begin{equation*}
    \underset{n\rightarrow\infty}{\lim}\Sigma_{\sigma}(F,A_{n})=F(0).
\end{equation*}

Whence
\begin{equation*}
    \{A_{n}\}\sim_{\sigma}(0,G),
\end{equation*}
for any domain $G$ satisfying the requirements of Definition $\ref{def-distribution}$.
\end{example}

\begin{example} Consider a 3-level $\alpha$-Toeplitz matrix $A$ where $\alpha=(\alpha_{1},\alpha_{2},\alpha_{3})=(0,1,2)$
\begin{eqnarray*}
  A&=&\left[\left[\left[a_{(r_{1}-0\cdot s_{1},r_{2}-1\cdot s_{2},r_{3}-2\cdot s_{3})}
      \right]_{r_{3},s_{3}=0}^{n_{3}-1}\right]_{r_{2},s_{2}=0}^{n_{2}-1}
      \right]_{r_{1},s_{1}=0}^{n_{1}-1}\\
   &=&\left[\left[\left[a_{(r_{1},r_{2}-s_{2},r_{3}-2s_{3})}
      \right]_{r_{3},s_{3}=0}^{n_{3}-1}\right]_{r_{2},s_{2}=0}^{n_{2}-1}
      \right]_{r_{1}=0}^{n_{1}-1}.
\end{eqnarray*}

The procedure is the same as in the previous example of an $\alpha$-circulant matrix, but in this case we do not need to permute the vector $\alpha$ since the only component equal to zero is already in first position. For $r_{1}=0,1,...,n_{1}-1$, let us set
\begin{eqnarray*}
    T_{r_{1}}=\left[\left[a_{(r_{1},r_{2}-s_{2},r_{3}-2s_{3})}
      \right]_{r_{3},s_{3}=0}^{n_{3}-1}\right]_{r_{2},s_{2}=0}^{n_{2}-1},
\end{eqnarray*}
then $T_{r_{1}}$ is a 2-level $\alpha^+$-Toeplitz matrix with $\alpha^+=(1,2)$ and of partial
sizes $n[>0]=(n_{2},n_{3})$ and
\begin{eqnarray*}
  A=\left[\begin{array}{cccc}
       T_{0} & T_{0} & \cdots & T_{0} \\
       T_{1} & T_{1} & \cdots & T_{1} \\
       \vdots & \vdots & \vdots & \vdots \\
       T_{n_{1}-1} & T_{n_{1}-1} & \cdots & T_{n_{1}-1}
    \end{array}\right].
\end{eqnarray*}
The latter is a block matrix with constant blocks on each row.
From formula $(\ref{1})$, the singular values of $A$ are the square root of the eigenvalues of $A^{*}A$:
\begin{eqnarray*}
  A^{*}A &=&\left[\begin{array}{cccc}
                T_{0}^{\ast} & T_{1}^{\ast} & \cdots & T_{n_{1}-1}^{\ast} \\
                T_{0}^{\ast} & T_{1}^{\ast} & \cdots & T_{n_{1}-1}^{\ast} \\
                \vdots & \vdots & \vdots & \vdots \\
                T_{0}^{\ast} & T_{1}^{\ast} & \cdots & T_{n_{1}-1}^{\ast} \\
                \end{array}\right]\left[\begin{array}{cccc}
                T_{0} & T_{0} & \cdots & T_{0} \\
                T_{1} & T_{1} & \cdots & T_{1} \\
                \vdots & \vdots & \vdots & \vdots \\
                T_{n_{1}-1} & T_{n_{1}-1} & \cdots & T_{n_{1}-1}
                \end{array}\right] \\
    &=&\left[\begin{array}{cccc}
          \overset{n_{1}-1}{\underset{j=0}\sum}T_{j}^{\ast}T_{j} &
          \overset{n_{1}-1}{\underset{j=0}\sum}T_{j}^{\ast}T_{j} & \cdots &
          \overset{n_{1}-1}{\underset{j=0}\sum}T_{j}^{\ast}T_{j} \\
          \overset{n_{1}-1}{\underset{j=0}\sum}T_{j}^{\ast}T_{j} &
          \overset{n_{1}-1}{\underset{j=0}\sum}T_{j}^{\ast}T_{j} & \cdots &
          \overset{n_{1}-1}{\underset{j=0}\sum}T_{j}^{\ast}T_{j} \\
          \vdots & \vdots & \vdots & \vdots\\
          \overset{n_{1}-1}{\underset{j=0}\sum}T_{j}^{\ast}T_{j} &
          \overset{n_{1}-1}{\underset{j=0}\sum}T_{j}^{\ast}T_{j} & \cdots &
          \overset{n_{1}-1}{\underset{j=0}\sum}T_{j}^{\ast}T_{j} \\
        \end{array}\right] \\
   &=&  \underbrace{\left[\begin{array}{cccc}
         1 & 1 & \cdots& 1 \\
         1 & 1 & \cdots& 1 \\
         \vdots & \vdots & \vdots & \vdots \\
         1 & 1 &\cdots& 1 \\
       \end{array}\right]}_{\textrm{$n_{1}$ times}}\otimes\overset{n_{1}-1}{\underset{j=0}\sum}T_{j}^{\ast}T_{j}\\
   &=& J_{n_{1}}\otimes\overset{n_{1}-1}{\underset{j=0}\sum}T_{j}^{\ast}T_{j}.
\end{eqnarray*}

Therefore
\begin{equation}\label{AstarAhat}
    {\rm Eig}(A^{*}A)={\rm Eig}\left(J_{n_{1}}\otimes\overset{n_{1}-1}{\underset{j=0}\sum}T_{j}^{\ast}T_{j}\right),
\end{equation}
where
\begin{eqnarray}\label{eigJ1}
  {\rm Eig}(J_{n_{1}}) = \{0,n_{1}\},
\end{eqnarray}
because $J_{n_{1}}$ is a
matrix of rank 1, so it has all eigenvalues equal to zero except one
eigenvalue equal to ${\rm tr}(J_{n_{1}})=n_{1}$ (${\rm tr}$ is the trace of a matrix). If  we put
\begin{eqnarray*}
  \lambda_{k}=\lambda_{k}\left(\overset{n_{1}-1}{\underset{j=0}\sum}T_{j}^{\ast}T_{j}\right),\qquad k=0,\ldots,n_{3}n_{2}-1,
\end{eqnarray*}
by exploiting basic properties of the tensor product and taking into
consideration $(\ref{AstarAhat})$ and $(\ref{eigJ1})$ we find
\begin{eqnarray}\label{eigAAhat1}
  \lambda_{k}(A^{*}A)&=&n_{1}\lambda_{k},\qquad k=0,\ldots,n_{3}n_{2}-1,\\
  \label{eigAAhat2}\lambda_{k}(A^{*}A)&=&0,\qquad k=n_{3}n_{2},\ldots,n_{3}n_{2}n_{1}-1.
\end{eqnarray}

From $(\ref{eigAAhat1})$, $(\ref{eigAAhat2})$ and $(\ref{1}),$ one obtains that
the singular values of $A$ are given by
\begin{eqnarray}
  \sigma_{k}(A)&=&\sqrt{n_{1}\lambda_{k}},\qquad k=0,\ldots,n_{3}n_{2}-1,\\
  \sigma_{k}(A)&=&0,\qquad k=n_{3}n_{2},\ldots,n_{3}n_{2}n_{1}-1.
\end{eqnarray}
and, since $\overset{n_{1}-1}{\underset{j=0}\sum}T_{j}^{\ast}T_{j}$ is a positive semidefinite matrix, from $(\ref{bmezzi})$ we can write
\begin{eqnarray*}
  \sigma_{k}(A)&=&\sqrt{n_{1}}\widetilde{\sigma}_{k},\qquad k=0,\ldots,n_{3}n_{2}-1,\\
  \sigma_{k}(A)&=&0,\qquad k=n_{3}n_{2},\ldots,n_{3}n_{2}n_{1}-1,
\end{eqnarray*}
where $\widetilde{\sigma}_{k}$ denotes the generic singular value of $\left(\overset{n_{1}-1}{\underset{j=0}\sum}T_{j}^{\ast}T_{j}\right)^{1/2}$.

Regarding the distribution in the sense of singular values, by invoking exactly the same argument as in the above example for $\alpha$-circulant matrix, we deduce that
\begin{equation*}
    \{A_{n}\}\sim_{\sigma}(0,G),
\end{equation*}
for any domain $G$ satisfying the requirements of Definition $\ref{def-distribution}$.
\end{example}

\begin{example} Let us see what happens when the vector $\alpha$ has only one component different from zero. Let $n=(n_{1},n_{2},\ldots,n_{d})$ and $\alpha=(0,\ldots,0,\alpha_{k},0,\ldots,0)$, $\alpha_{k}>0$; in this case we can give an explicit formula for the singular values of the $d$-level $\alpha$-circulant matrix. For convenience and without loss of generality we take $\alpha=(0,\ldots,0,\alpha_{d})$ (with all zero components in top positions, otherwise we use a permutation). From  $\ref{nonnegative-vs-positive}$, the singular values of
$A_{n}=[a_{(r-\alpha\circ s)\textrm{ mod $n$}}]_{r,s=\underline{0}}^{n-e}$
are zero except for few of them given by $\sqrt{\hat n[0]} \sigma$ where, in our case, $\hat n[0]=n_{1}n_{2}\cdots n_{d-1}$, $n[0]=(n_{1},n_{2},\ldots,n_{d-1})$, and $\sigma$ is any singular value of the matrix
\begin{equation*}
\left(\sum_{j=\underline{0}}^{n[0]-e}C_j^*C_j\right)^{1/2},
\end{equation*}
where $C_j$ is an $\alpha_{d}$-circulant matrix of dimension $n_{d}\times n_{d}$ whose expression is
\begin{eqnarray*}
C_j=\left[a_{(r-\alpha \circ s)\ {\rm mod}\,
n}\right]_{r_{d},s_{d}=0}^{n_{d}-1}&=&\left[a_{(r_{1},r_{2},\ldots,r_{d-1},(r_{d}-\alpha_{d}s_{d})\ {\rm mod}\,
n_{d})}\right]_{r_{d},s_{d}=0}^{n_{d}-1}\\
&=&\left[a_{(j,(r_{d}-\alpha_{d}s_{d})\ {\rm mod}\,
n_{d})}\right]_{r_{d},s_{d}=0}^{n_{d}-1},
\end{eqnarray*}
with $(r_{1},r_{2},\ldots,r_{d-1})=j$. For $j=\underline{0},\ldots,n[0]-e$, if $C_{n_{d}}^{(j)}$ is the circulant matrix which has as its first column the vector $a^{(j)}=[a_{(j,0)},a_{(j,1)},\ldots,a_{(j,n_{d}-1)}]^{T}$ (which is the first column of the matrix $C_{j}$), $C_{n_{d}}^{(j)}=[a_{(j,(r-s)\textrm{ mod $n_{d}$})}]_{r,s=0}^{n_{d}-1}=F_{n_{d}}D_{n_{d}}^{(j)}F_{n_{d}}^{*}$, with $D_{n_{d}}^{(j)}=diag(\sqrt{n_{d}}F_{n_{d}}^{*}a^{(j)})$, then, from $(\ref{CastC})$, $(\ref{0})$, and $(\ref{Vi})$,
it is immediate to verify that
\begin{eqnarray*}
\sum_{j=\underline{0}}^{n[0]-e}C_j^*C_j&=&\sum_{j=\underline{0}}^{n[0]-e}(F_{n_{d}}D_{n_{d}}^{(j)}F_{n_{d}}^{*}Z_{n_{d},\alpha_{d}})^{*}
(F_{n_{d}}D_{n_{d}}^{(j)}F_{n_{d}}^{*}Z_{n_{d},\alpha_{d}})\\
&=&\sum_{j=\underline{0}}^{n[0]-e}(F_{n_{d}}^{*}Z_{n_{d},\alpha_{d}})^{*}(D_{n_{d}}^{(j)})^{*}D_{n_{d}}^{(j)}(F_{n_{d}}^{*}Z_{n_{d},\alpha_{d}})\\
&=&(F_{n_{d}}^{*}Z_{n_{d},\alpha_{d}})^{*}\left(\sum_{j=\underline{0}}^{n[0]-e}(D_{n_{d}}^{(j)})^{*}D_{n_{d}}^{(j)}\right)(F_{n_{d}}^{*}Z_{n_{d},\alpha_{d}}).
\end{eqnarray*}

Now, if we put $n_{d,\alpha}=\frac{n_{d}}{(n_{d},\alpha_{d})}$ and
\begin{eqnarray*}
   q_{s}^{(j)}&=&|D_{n_{d}}^{(j)}|_{s,s}^{2}=(D_{n_{d}}^{(j)})_{s,s}\cdot\overline{(D_{n_{d}}^{(j)})_{s,s}},
   \quad s=0,1,\ldots,n_{d}-1,\\
   \Delta_{l}&=&\left[\begin{array}{cccc}
       \overset{n[0]-e}{\underset{j=\underline{0}}\sum} q_{(l-1)n_{d,\alpha}}^{(j)} & & & \\
   & \overset{n[0]-e}{\underset{j=\underline{0}}\sum} q_{(l-1)n_{d,\alpha}+1}^{(j)} & & \\
   & & \ddots & \\
   & & & \overset{n[0]-e}{\underset{j=\underline{0}}\sum} q_{(l-1)n_{d,\alpha}+n_{d,\alpha}-1}^{(j)}
   \end{array}\right]\in \mathbb{C}^{n_{d,\alpha}\times
   n_{d,\alpha}},\\
\end{eqnarray*}
for $l=1,2,\ldots,(n_{d},\alpha_{d})$, then, following the same reasoning employed for proving formula $(\ref{eigg})$, we infer
\begin{equation*}
    {\rm Eig}\left(\sum_{j=\underline{0}}^{n[0]-e}C_j^*C_j\right)=\frac{1}{(n_{d},\alpha_{d})}{\rm Eig}\left(J_{(n_{d},\alpha_{d})}
    \otimes \overset{(n_{d},\alpha_{d})}{\underset{l=1}\sum}\Delta_{l}\right),
\end{equation*}
where
\begin{eqnarray*}
  J_{(n_{d},\alpha_{d})}&=&\underbrace{\left[\begin{array}{cccc}
   1 & 1 & \cdots & 1\\
   1 & 1 & \cdots & 1\\
   \vdots & \vdots & \vdots & \vdots\\
   1 & 1 & \cdots & 1
\end{array}\right]}_{\textrm{$(n_{d},\alpha_{d})$ times}},\\
  \frac{1}{(n_{d},\alpha_{d})}{\rm Eig}(J_{(n_{d},\alpha_{d})}) &=& \{0,1\},
\end{eqnarray*}
and
\begin{eqnarray*}
    \overset{(n_{d},\alpha_{d})}{\underset{l=1}\sum}\Delta_{l}&=&
    \overset{(n_{d},\alpha_{d})}{\underset{l=1}\sum}{\rm diag}\left(\overset{n[0]-e}{\underset{j=\underline{0}}\sum}q_{(l-1)n_{d,\alpha}+k}^{(j)};\text{\,\,}k=0,1,\ldots,n_{d,\alpha}-1\right)\\
    &=&{\rm diag}\left(\overset{(n_{d},\alpha_{d})}{\underset{l=1}\sum}
    \overset{n[0]-e}{\underset{j=\underline{0}}\sum}q_{(l-1)n_{d,\alpha}+k}^{(j)};\text{\,\,}
    k=0,1,\ldots,n_{d,\alpha}-1\right).
\end{eqnarray*}

Consequently, since $\overset{(n_{d},\alpha_{d})}{\underset{l=1}\sum}\Delta_{l}$ is a diagonal
matrix, and by exploiting basic properties of the tensor product, we find
\begin{eqnarray*}
  \lambda_{k}\left(\sum_{j=\underline{0}}^{n[0]-e}C_j^*C_j\right)&=&\overset{(n_{d},\alpha_{d})}
  {\underset{l=1}\sum}\overset{n[0]-e}{\underset{j=\underline{0}}\sum}q_{(l-1)n_{d,\alpha}+k}^{(j)}, \quad k=0,1,\ldots,n_{d,\alpha}-1,\\
  \lambda_{k}\left(\sum_{j=\underline{0}}^{n[0]-e}C_j^*C_j\right)&=& 0,\qquad k=n_{d,\alpha},\ldots,n_{d}-1.
\end{eqnarray*}

Now, since $\sum_{j=\underline{0}}^{n[0]-e}C_j^*C_j$ is a positive semidefinite matrix, from $(\ref{bmezzi})$ we finally have
\begin{eqnarray*}
  \sigma_{k}\left(\left(\sum_{j=\underline{0}}^{n[0]-e}C_j^*C_j\right)^{1/2}\right)&=&\sqrt{\overset{(n_{d},\alpha_{d})}
  {\underset{l=1}\sum}\overset{n[0]-e}{\underset{j=\underline{0}}\sum}q_{(l-1)n_{d,\alpha}+k}^{(j)}}, \quad k=0,1,\ldots,n_{d,\alpha}-1,\\
  \sigma_{k}\left(\left(\sum_{j=\underline{0}}^{n[0]-e}C_j^*C_j\right)^{1/2}\right)&=& 0,\qquad k=n_{d,\alpha},\ldots,n_{d}-1.
\end{eqnarray*}
\end{example}

\section{Conclusions and future work}\label{sec:fin}

In this paper we have studied in detail the singular values of
$\alpha$-circulant matrices and we have identified the joint
asymptotic distribution of $\alpha$-Toeplitz sequences associated
with a given integrable symbol. The generalization to the multilevel
block setting has been sketched together with some intriguing
relationships with the design of multigrid procedures for structured
linear systems. The latter point deserves more attention and will be
the subject of future researches. We also would like to study the
more involved eigenvalue/eigenvector behavior both for
$\alpha$-circulant and $\alpha$-Toeplitz structures.

\end{document}